\documentclass[3p,times]{elsarticle}

\usepackage[utf8]{inputenc}

\usepackage{graphicx}
\usepackage[dvipsnames]{xcolor}

\usepackage{mathrsfs}
\usepackage{amsfonts}
\usepackage{amssymb}
\usepackage{amsmath}
\usepackage{amsthm}
\usepackage{dsfont}
\usepackage{rotating}
\usepackage{multirow}
\usepackage{subfigure}
\usepackage{comment}
\usepackage{cancel}
\usepackage[noabbrev,capitalise]{cleveref}

\usepackage{mathtools}
\usepackage{pgffor}
\usepackage{tikz}
\usepackage{pgfplots}
\pgfplotsset{compat = newest}
\usepackage{booktabs}
\usepackage{rotating}
\usepackage{lscape}
\usepackage{graphics}

\DeclareMathAlphabet\mathbfcal{OMS}{cmsy}{b}{n}  

\newcommand{\x}{\mathbf{x}}
\newcommand{\q}{\mathbf{q}}
\newcommand{\CFL}{\textnormal{CFL}}

\newcommand{\R}{\mathds{R}}
\newcommand{\f}{\mathbf{f}} 

\newcommand{\Q}{\mathbf{Q}} 
 
\newcommand{\E}{\mathcal{E}}
\newcommand{\F}{\mathbf{F}}
\newcommand{\B}{\mathbf{B}}

\newcommand{\N}{\mathbb{N}}
\newcommand{\n}{\mathbf{n}}

\renewcommand{\u}{\mathbf{u}}

\renewcommand{\vec}[1]{\mathbf{#1}}

\newtheorem{remark}{Remark}[section]

\newtheorem{proposition}{Proposition}[section]

\newcommand{\diff}[1]{{\mathrm{d}{#1}}}

\newcommand{\partialfra}[1]{\frac{\partial}{\partial x^{#1}}}
\newcommand{\minmod}{\text{minmod}}

\newcommand{\RIcolor}[1]{{\leavevmode\color{black} #1}}
\newcommand{\RIIcolor}[1]{{\leavevmode\color{black} #1}}
\newcommand{\RIIIcolor}[1]{{\color{black} #1}}

\begin{document}
	
\begin{frontmatter}
	
	\journal{Applied Mathematics and Computation} 
	
	\title{Well balanced finite volume schemes for shallow water equations on manifolds } 

	\cortext[mycorrespondingauthor]{Corresponding author}
	
	\author[Inria]{Michele Giuliano Carlino\corref{mycorrespondingauthor}}
	\ead{michele-giuliano.carlino@inria.fr }
	\author[Inria]{Elena Gaburro}
	\ead{elena.gaburro@inria.fr}

	\address[Inria]{Inria, Univ. Bordeaux, CNRS, Bordeaux INP, IMB, UMR 5251, 200 Avenue de la Vieille Tour, 33405 Talence cedex, France }

	%
	\begin{abstract}
		
		In this paper we propose a novel second-order accurate well balanced scheme for shallow water equations in general covariant coordinates over manifolds. 
		In our approach, once the gravitational field is defined for the specific case, one equipotential surface is detected and parametrized by a frame of general covariant coordinates. This surface is the manifold whose covariant parametrization induces a metric tensor. The model is then re-written in a hyperbolic form with a tuple of conserved variables composed both of the evolving physical quantities and the metric coefficients. 
		This formulation allows the numerical scheme to automatically compute the curvature of the manifold as long as the physical variables are evolved. 
		
		In a classical well balanced formulation, the knowledge of a given equilibrium is exploited for evolving the specific state seen as the sum of the equilibrium and the fluctuations of the state around the equilibrium itself. 
		On the contrary, the numerical approach proposed here is \textit{automatically well balanced} for the \textit{water at rest} solution (defined by zero velocity and constant free surface) for \textit{general manifolds}
		without having to exploit the exact equilibrium profile during the computation.		
		In particular, this numerical strategy allows to preserve the accuracy of the water at rest equilibrium at machine precision and on large timescales, even for non-smooth bottom topographies. As a matter of fact, the proposed local polynomial reconstruction, needed at each time step by the scheme for evolving the state, is built in order to automatically cancel any numerical error committed for describing jumps in the bathymetry. As a further effect, also out of the equilibrium, the typical spurious non-physical oscillations of the recovered numerical solution in a neighborhood of the discontinuities are healed. Thus, once the information on both the bathymetry and the metric of the manifold are properly collected in the flux terms and in the nonconservative products, the resulting numerical strategy turns out to be more accurate in finding also non-equilibrium solutions. 
		
		Numerical results close the work. After having proved that the scheme is second-order accurate, we test the recovery of the water at rest equilibrium at machine precision for different bottom topographies (eventually discontinuous) with various metric on very long simulation times. In particular, after having considered smooth bathymetries with compact support, we relax the assumption on the continuity of water at rest equilibrium with discontinuous bathymetries disturbed by a white noise.		
		  
	\end{abstract}	
	%
\end{frontmatter}


\section{Introduction} \label{sec.intro}

The goal of this article consists in developing a second-order accurate scheme for Shallow Water (SW) equations  in General Covariant Coordinates (GCCs) over a manifold. The proposed scheme is able to preserve the special equilibrium of water at rest (i.e. zero velocity and constant free surface) over long times. The considered physical model, discretized by the scheme, takes into account the particular shape of the gravitational field leading the dynamics of the movement of the fluid of the test case in consideration.

SW equations have been a useful tool for describing various phenomena related to fluid geophysics, from the study of flows in a lake to the movement of water in a river~\cite{lanzoni2006long, zolezzi2001downstream, fraccarollo2002riemann, george2014depth, iverson2014depth}, to appropriately describing meteorological~\cite{holton1973introduction, jacobs1997diurnal} or oceanic~\cite{higdon2006numerical, cotter2009mixed} phenomena over long distances or providing a useful risk analysis tool for avalanches~\cite{gray1999gravity} or tsunamis~\cite{segur2007waves}. Although all of these physical phenomena respond to the same type of mathematical description, the input datum of the gravitational field can change with respect to the spatial scale and with respect to the bottom topography under consideration~\cite{staniforth2015dynamically, benard2015assessment, staniforth2015consistent}. 
Indeed, starting with the shallow water assumption for which fluid particles have negligible motion along the gravitational field lines~\cite{pedlosky1987geophysical} (shallow water hypothesis), consider an immersed manifold that is a significant surface perpendicular to the gravitational field lines. From the manifold, a certain metric is deduced that allows us to describe the differential operator of SW equations. 
Such an operator will then be strongly dependent on the metric, thus on the curvature of the manifold, and will change from case to case. Concrete examples may be those involving the fluid dynamic study of a lake, in which gravity can be modeled as perpendicular to a horizontal plane, or oceanic dynamics that take into account a spherical distribution of gravity~\cite{williamson1992standard, kolar1994shallow, arpaia2020well, baldauf2020discontinuous} up to more detailed descriptions that, accounting for the different local mass placement of the planet, consider gravity as a field that varies in modulus and direction with respect to the physical-geographical characteristics of the Earth's geoid~\cite{staniforth2015shallow, staniforth2015geophysically}. 

The starting mathematical model of this paper thus considers as input the metric description of the manifold, and once this is defined as initial datum, it is appropriately integrated into the differential operator leading to a hyperbolic differential system that evolves the solution with respect to the geometric characteristics of the manifold itself in an autonomous way.

In general, for small variations in the bathymetry, also for large distances (e.g. oceanic distances), the fluid depth is evaluated in vertical direction, which does not necessarily coincide with the direction of the gravitational forces. This violates the previously stated shallow water hypothesis. 
Although such an approximation might seem reasonable with respect to the variation of the gravitational field perpendicular to a certain surface, Bachini and Putti in~\cite{bachini2020geometrically} have shown that the SW equations resulting i) from the integration of the Navier-Stokes equations along gravitational field lines and ii) along the vertical reference line lead to significantly different numerical results for the same physical phenomenon, even for bathymetries that vary very slightly in space. Concerning spherical gravitational field, similar differences have been recorded also by Chow in~\cite{chow1959open}. Therefore, for a numerical output to be as accurate as possible, the shape of the gravitational field must be properly accounted for in the mathematical model describing a certain phenomenon.

In addition to the SW classical model, in Cartesian coordinates, 
in literature it is possible to find SW models that are described by only one fixed frame of general covariant coordinates, for example see~\cite{castro2017well, arpaia2020well} for spherical coordinates. 
A particular case is given by Baldauf in~\cite{baldauf2020discontinuous}, from which our work starts (see Sec.~\ref{sec.model}). In this last cited article, despite a general formulation in a general covariant frame is provided for describing the shallow water system, all the geometrical features of the accounted spherical and elliptical manifolds are considered in the system as input data, included the curvatures. 
Instead, the approach we propose here is different, because the hyperbolic system we formulate, in addition to the physical quantities, only needs the metric to be given in input (and the metric will be part of the set of conserved variables). 
Then, once the metric is defined, the system automatically computes all the other geometric terms describing the curvature of the manifold (for example the Christoffel symbols). 
So, at the same time, our model evolves the conserved physical variables and computes the specific geometric features of the manifold, being thus more general and flexible.

This approach is strongly inspired by the work by Gaburro, Dumbser \& Castro in~\cite{gaburro2021well}. Although this paper deals with numerical schemes on general relativity, basically the gravitational terms in fixed backgrounds are not considered in algebraic form as a source term but, through the use of nonconservative terms, they are expressed as a function of the spatial derivatives of appropriate conserved variables. Equivalently in our context, the terms related to the spatial derivatives in general covariant coordinates of the covariant formulation of the metric
define the nonconservative terms of the hyperbolic system. 
Thus, the \textit{tuple of conserved variable} is composed not only of the \textit{physical variables} of the system but also of the \textit{geometric description of the manifold}. The advantage of this formulation is that, therefore, the spatial curvature is calculated discretely in a system-autonomous form and taken into account in the evolution of physics.

This is made possible by considering a general covariant formulation of the problem. 
In the literature also another approach via Local Covariant Coordinates (LCC) can be found (e.g. see~\cite{baldauf2020discontinuous, fent2018modeling, arpaia2020well}). 
One advantage of LCC strategy is that eventual critical points for the metric (i.e. those points of the computational domain over which the determinant of the metric vanishes, e.g. the poles onto a sphere in classical longitudinal-latitudinal spherical mapping) are avoided because any cell (also those containing these pathological points) are mapped into a reference system where the numerical solution is recovered. However, the numerical scheme that is constructed increases in complexity because this reference map is added to the one that parameterizes the manifold in the computational domain. 
So, specific strategies need to be employed to recover the possible continuity of the solution or alternatively the conservativity of the numerical fluxes. Moreover, in a local covariant coordinates approach it becomes more challenging to deal with discontinuous bathymetry because the local map along the jumps usually becomes singular unless special \textit{ad hoc} constructed numerical manipulations are employed~\cite{valiani2017momentum}.

\medskip 

Besides, the objective of preserving over long simulation times one or more equilibria of the model, 
can be assessed via the well balancing feature of the numerical method. 
When a scheme is well balanced, it eliminates spurious numerical modes, it reduces the numerical dissipation and small-scale phenomena may be more evident on larger scales. 
In~\cite{bermudez1994upwind}, Bermudez and Vazquez use shallow water system to well balance it with respect to water at rest equilibrium in a rigorous mathematical fashion. Afterwards, concerning this topic, the literature has expanded; as examples we cite~\cite{leveque1998balancing, perthame2001kinetic, gosse2001well, rebollo2003family, bouchut2004nonlinear, audusse2004fast, tang2004gas, castro2008well, kappeli2014well, chandrashekar2015second, castrodeluna2017well, gaburro2017direct, gaburro2018well, klingenberg2019arbitrary, guerrero2020second, thomann2020all, berberich2020high, gomez2021high, guerrero2021well, berberich2021high}. A comprehensive literature review on the topic is given by Castro \& Par{\'e}s in~\cite{castro2020well}. With respect to well balancing for SW, in addition,~\cite{lukavcova2007well, noelle2006well,noelle2007high,russo2008high,kazolea2013well,xing2006high,abgrall2022hyperbolic} can be cited.

A crucial aspect of this article is to relax the continuity hypothesis of the solution. Typically, equilibrium is considered to be a smooth steady solution. But in general this assumption cannot be considered. In practical situations it is not always possible to detect the bathymetry in continuous form. Therefore, even where the bathymetry is at least Lipshitz-continuous, it is discretized as discontinuous. For that reason, the depth of the fluid would also be discontinuous. In this paper we propose a way such that any numerical error due to jumps in bathymetry automatically cancels when the solution evolves. This approach not only preserves the water at rest equilibrium but also allows us to recover a numerical solution out the equilibrium that is not affected by any bad approximations of jumps in the bottom.

\medskip
	
\RIIIcolor{
The paper is organized as follows. The original shallow water model on manifolds in general covariant coordinates is presented in Sec.~\ref{sec.model}. In this section the hyperbolic system takes in input the metric as well as all the geometrical features of the manifold. After some physical remarks in order to properly define the meaning of the manifold in this context, in Sec.~\ref{subsec.WRcontinuum} the water at rest equilibrium is introduced with respect to a general covariant frame. Afterwards, Sec.~\ref{sec.ourmodel} is devoted to the derivation of the final formulation of the novel shallow water governing model in general covariant coordinates only depending on the metric. In particular, the new hyperbolic system, formally presented in Sec.~\ref{subsec.formulation}, is able to evolve the physical variables and, at the same time, to automatically compute the geometrical features of the manifold. The adopted numerical scheme is described in Sec. \ref{sec.method} and its well balanced formulation is analysed in Sec. \ref{sec.WB}. Numerical validations are in Sec.~\ref{sec.test}. The paper is closed by conclusions and future perspectives in Sec.~\ref{sec.concl}.
}

\section{Shallow water equations on manifolds in general covariant coordinates} \label{sec.model}

\begin{figure}
	\centering
	\begin{tikzpicture} 
		\begin{axis}[
			hide axis,
			enlargelimits={abs=10pt},
			view = {50}{40}
			]

			\addplot3 [domain=0:10, y domain = 0:10,
          		quiver = {
          			u={-2*(x)/sqrt(4*x^2 + 4*y^2 +1)},
          			v={-2*(y)/sqrt(4*x^2 + 4*y^2 +1)},
          			w={-2000},
          			scale arrows=.1
          		},
          		-stealth,
          		samples=7,
          		red] ({x}, {y}, {(-(x^2 + y^2) - 150)});	
          	\draw (9, 5, -600) node[right] {\textcolor{red}{$\mathbf{g}$}}	;	
			
			\addplot3 [
    				domain=-10:10,
    				domain y = -10:10,
    				only marks,
    				mark size=0.5pt,
    				gray] {-(x^2 + y^2) - 150};
    			\draw (10, 10, -370) node[right] {\textcolor{gray}{$\mathcal{M}$}}	;
    				
    				\addplot3 [
    				domain=-10:10,
    				domain y = -10:10,
    				only marks,
    				mark size=0.5pt,
    				mesh,
    				draw = black] {-(x^2 + y^2) + 15*x*(y > 0) + 16*y*(x < 1.5) + 150};
    				\draw (9, 9, 0) node[right] {$b$}	;
    				
    				\addplot3 [
    				domain=-10:10,
    				domain y = -10:10,
    				surf, 
    				color=blue, 
    				opacity=0.4,
    				faceted color=blue!40] {-(x^2 + y^2) + 50*cos(x*y*5) + 450};
    				\addplot3 [domain=-10:10, y domain = 0:0, color=blue, samples=100] (x,-10,{-(x^2 + 100) + 50*cos(-x*50) + 450});
          		\addplot3 [domain=-10:10, y domain = 0:0, color=blue, samples=100] (x,10,{-(x^2 + 100) + 50*cos(x*50) + 450});
          		\addplot3 [domain=-10:10, y domain = 0:0, color=blue, samples=100] (-10,x,{-(x^2 + 100) + 50*cos(-x*50) + 450});
          		\addplot3 [domain=-10:10, y domain = 0:0, color=blue, samples=100] (10,x,{-(x^2 + 100) + 50*cos(x*50) + 450});
          		\draw (9, 10, 300) node[above] {\textcolor{blue}{$\eta$}}	;

          		\draw[->, red] (-9, -9, 600) -- (-8.5, -8.5, -600);
          		\draw plot [mark=*, mark size=1] coordinates{(-8.8917, -8.8917, 340)}; 
          		\draw plot [mark=*, mark size=1] coordinates{(-8.6792, -8.6792, -170)}; 
          		\draw[densely dotted, black, thick] (-8.8917, -8.8917, 340) -- (-8.6792, -8.6792, -170);
          		\draw[black, thick, ->] (-8.8917, -8.8917, 340) -- (-8.8917-1.8, -8.8917, 340);
				\draw[black, thick, ->] (-8.8917, -8.8917, 340) -- (-8.8917, -8.8917+3, 340);	
				\draw (-8.5, -8.5, -600) node[below] {\textcolor{red}{$\mathbf{g}$}}	;
				\draw (-9.5, -9, 100) node[below] {$h$}	;
				\draw (-8.8917-1.8, -8.8917, 340) node[above] {$u^1$}	;
				\draw (-8.8917, -8.8917+3, 340) node[below] {$u^2$}	;

				\draw[<-, thick] plot [smooth, domain=0:1, samples=10] ( {8+5*\x}, {-12}, {-((8+5*\x)^2 + 110.25) - 150} ); 
				\draw[->, thick] plot [smooth, domain=0:1, samples=10] ( {12}, {-12.5+5*\x}, {-((-12.5+5*\x)^2 + 110.25) - 150});
				\draw (8, -12, -324) node[left] {$x^1$}	; 
				\draw (12, -7.5, -316.5) node[right] {$x^2$}	;

		\end{axis}
	\end{tikzpicture}
	\caption{Sketch of the geometric setting of SW in GCC~\eqref{eq.baldauf}. The gravitational field $\mathbf{g}$ defines an equipotential surface $\mathcal{M}$. Thus, this surface is the manifold over which the problem is solved and it is fully described by the frame of GCC $\x=(x^1, x^2)$. The bathymetry is given by the bottom surface $b$, eventually discontinuous (meshed black surface). The free surface $\eta$ of the fluid is the blue surface. Consequently, the fluid depth $h$ is the length of the gravitational field line conducted from the free surface $\eta$ up to the bottom $b$. The velocity $\u = (u^1, u^2)$ of the fluid is computed along the tangent plane of the manifold $\mathcal{M}$ and measured along the principal directions given by the GCC $\x$.}
	\label{fig.sketch_model}
\end{figure}
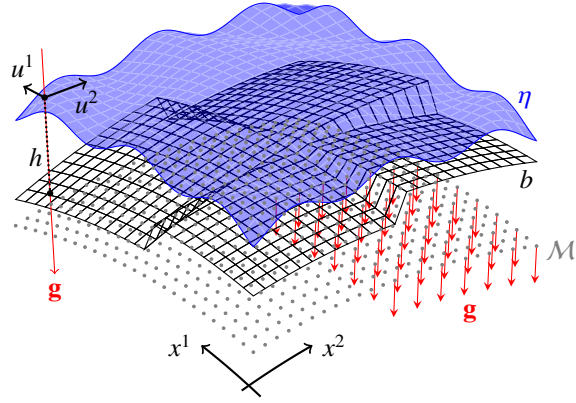

The considered differential system for shallow water equations on an arbitrary $d$-dimensional manifold $\mathcal{M}$ (with $d = 1,2$) in general covariant coordinates is originally introduced in~\cite{baldauf2020discontinuous} by Baldauf. For any nonnegative time $t$ and for any GCC $\x = (x^1, x^2)$ in the reference domain $\Omega \subset \R^d$, it reads
\begin{subequations} \label{eq.baldauf}
	\begin{alignat}{2}
		\frac{\partial h}{\partial t} + \nabla_j m^j &= 0, \label{eq.baldauf_mass} \\
		\frac{\partial m^i}{\partial t} + \nabla_j T^{ij} &= S^i, \label{eq.baldauf_mom}
	\end{alignat}
\end{subequations}
for $i,j = 1,2$. In the above system and in the sequel the Einstein notation for sum on repeated indexes (both superscript and subscript) is used, with exception if they are written in round brackets. In system~\eqref{eq.baldauf}, the evolving unknowns in time are the fluid depth $h$ (in the conservation of the mass~\eqref{eq.baldauf_mass}) and the mass fluxes $m^i = h u^i$ (in the balance of the momentum~\eqref{eq.baldauf_mom}), with $u^i$ representing the two velocity components of the fluid with respect to the principal directions defined by the generalized coordinates $x^i$, see Fig.~\ref{fig.sketch_model}. In the momentum equation in~\eqref{eq.baldauf}, the flux is defined by the stress tensor
\begin{equation} \label{eq.StressTesor}
	T^{ij} = \frac{m^i m^j}{h} + \frac{1}{2} g h^2 \gamma^{ij},
\end{equation}
composed of the sum between a kinetic and a pressure component. In the latter, the gravitational constant $g$ and the contravariant metric tensor $\gamma^{ij}$ intervene. The source term $S^i$ is defined as
\begin{equation} \label{eq.source}
	S^i = - g h \gamma^{ij} \partial_{x^j} b,
\end{equation}
where $b$ is the bathymetry. In the source term~\eqref{eq.source} any contribution from Coriolis, bottom-friction or viscous forces is neglected. The treatment of these further terms depends on the particular test case one wishes to study and is beyond the scope of this article. The reader is referred to~\cite{tort2014consistent, bresch2003existence, grundy1985approach, baldauf2020discontinuous, fent2018modeling} for more specific instances in which the contributions deduced from the physics of the case under consideration are added to the only gravitational contribution in~\eqref{eq.source}. \\
\indent Since the covariant derivative $\nabla_j$ expands in the sum of the derivatives with respect to the covariant coordinate $x^j$ and the additional metric terms that account for the curvature of the manifold via Christoffel symbols, equations~\eqref{eq.baldauf} can be equivalently rewritten as
\begin{subequations} \label{eq.baldauf1}
	\begin{alignat}{2}
		\frac{\partial h}{\partial t} + \frac{\partial m^j}{\partial x^j} + \Gamma^j_{jk} m^k &= 0, \label{eq.baldauf1_mass} \\
		\frac{\partial m^i}{\partial t} + \frac{\partial T^{ij}}{\partial x^j} + \Gamma^i_{jk} T^{kj} + \Gamma^j_{jk} T^{ik} &= S^i. \label{eq.baldauf1_mom}
	\end{alignat}
\end{subequations}
By definition, the Christoffel symbol $\Gamma^i_{jk}$ of second kind is
\begin{equation} \label{eq.christoffel}
	\Gamma^i_{jk} = \frac{1}{2} \gamma^{im} \bigg( \frac{\partial \gamma_{km}}{\partial x^j} + \frac{\partial \gamma_{jm}}{\partial x^k}- \frac{\partial \gamma_{jk}}{\partial x^m} \bigg),
\end{equation}
where $\gamma_{ij}$ is the covariant representation of the metric tensor. Thus, it is the inverse of the contravariant metric $\gamma^{ij}$, i.e. $\gamma^{ij} \gamma_{jk} = \delta^i_k$, with $\delta^i_k$ the Kronecker symbol.\\
\indent Due to the definition of the gravitational contribution in the source term~\eqref{eq.source}, the manifold over which system~\eqref{eq.baldauf1} is solved has to be an equipotential surface~\cite{baldauf2020discontinuous}. Extension to non equipotential manifolds for SW in GCC~\eqref{eq.baldauf} are discussed by Staniforth in~\cite{staniforth2015dynamically}. As a direct consequence, the gravitational field $\vec{g}$ (whose modulus is taken constantly equal to $g$ in this paper) is perpendicular to the manifold (because it is normal to any equipotential surface). Thus, the height $h$ is defined as the length of the gravitational field line conducted from the free surface $\eta = h + b$ of the fluid to the bottom, as it is sketched in Fig.~\ref{fig.sketch_model}.

\begin{remark}
	In this paper, the fluid depth $h$, velocity component $u^i$ (or the relative mass flux $m^i$), for $i = 1,2$, and the bathymetry $b$ are said \textit{physical} quantities; all terms involving the metric tensor (in contravariant or covariant form) are called \textit{metric} quantities.
\end{remark}

\subsection{The water at rest equilibria at continuum level} \label{subsec.WRcontinuum}
One goal of this article is to develop a well balanced numerical scheme for the shallow water equations in general covariant coordinates. In particular, the well-balancing of the method is built with respect to the special equilibrium of water at rest. Such an equilibrium is one possible solution in the uncountable set of solutions that satisfy the following proposition.
\begin{proposition}[Frame-free features of water at rest equilibria] \label{prop.wr}
	For problem~\eqref{eq.baldauf1} with a nonsingular metric $\gamma^{ij}$ and an assigned bathymetry $b$, the water at rest equilibrium consists in a fluid depth $h$ and a couple of velocity components $(u^1, u^2)$ such that
	\begin{subequations} \label{eq.WR}
		\begin{equation} \label{eq.WRfs}
			\eta = h + b \equiv \text{const. }
		\end{equation}
		\text{and}
		\begin{equation} \label{eq.WRvel}
			u^1 = u^2 \equiv 0,
		\end{equation}
	\end{subequations}
	for any positive time $t$ and for any $\x \in \Omega$.
\end{proposition}

\begin{proof}
	Since the velocity components $u^i$ vanish, mass fluxes $m^i$ are zero. Consequently, the mass conservation equation~\eqref{eq.baldauf1_mass} is automatically satisfied. Moreover, condition~\eqref{eq.WRvel} implies that the stress tensor~\eqref{eq.StressTesor} reduces just to its pressure component
	\begin{equation} \label{eq.StressTensor_pre}
		T^{ij} = \frac{1}{2} g h^2 \gamma^{ij}.
	\end{equation}
	By applying the chain rule for deriving~\eqref{eq.StressTensor_pre} w.r.t. to $x^j$, equations~\eqref{eq.baldauf1_mom} turn into
	\begin{equation} \label{eq.baldauf1_momWR}
		g \gamma^{ij} \frac{\partial h}{\partial x^j} + \frac{1}{2} g h^2 \gamma^{ij} + \Gamma^i_{jk} \bigg( \frac{1}{2} g h^2 \gamma^{kj} \bigg) + \Gamma^j_{jk} \bigg( \frac{1}{2} g h^2 \gamma^{ik} \bigg) = - g h \gamma^{ij} \frac{\partial b}{\partial x^j}.
	\end{equation}
	Condition~\eqref{eq.WRfs} implies that $\partial_{x^j}(h + b) = 0$, for $j = 1,2$. Consequently, equations~\eqref{eq.baldauf1_momWR} read
	\begin{equation} \label{eq.baldauf1_momWR1}
		\frac{1}{2} g h^2 \bigg( \frac{\partial \gamma^{ij}}{\partial x^j} + \Gamma^i_{jk} \gamma^{kj} + \Gamma^j_{jk} \gamma^{ik} \bigg) = 0, \quad \forall i \in \{1,2\}.
	\end{equation}
	The quantity into brackets in the left hand side of~\eqref{eq.baldauf1_momWR1} is the covariant derivative of the metric tensor $\gamma^{ij}$ w.r.t. to $x^j$, i.e. $\nabla_j \gamma^{ij}$. Since the connection is chosen to have zero covariant derivative of the metric (in order to allow parallel transport operations which preserve angles and lengths, see Sec. 3.1 of~\cite{wald2010general}), equations~\eqref{eq.baldauf1_momWR1} are verified. As a consequence, the momentum conservation equations~\eqref{eq.baldauf1_mom} are satisfied as well.\\
Finally, it is proved that equilibrium~\eqref{eq.WR} of water at rest is a solution for SW equations in GCC~\eqref{eq.baldauf1} for any nonsingular metric.   
\end{proof}

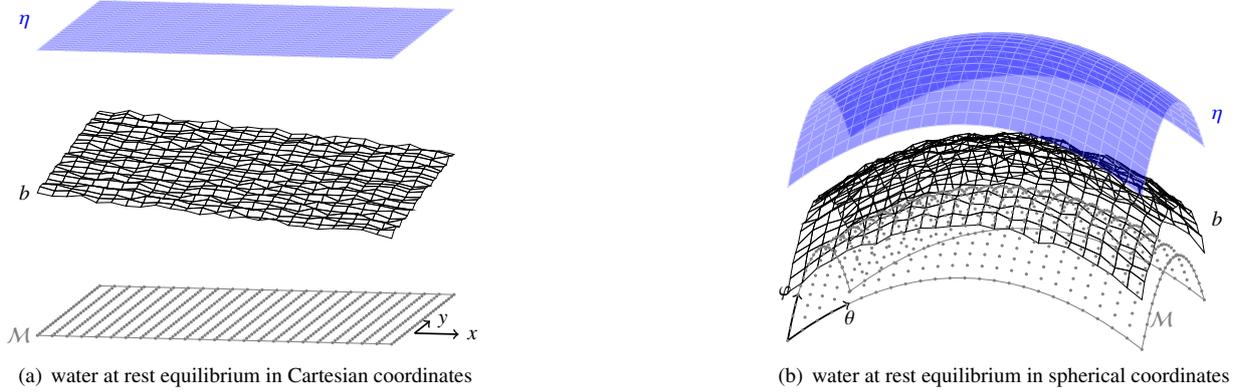
\begin{figure}
	\centering
	\subfigure[water at rest equilibrium in Cartesian coordinates]{
		\begin{tikzpicture}[scale=0.8]
			\begin{axis}[
			hide axis,
			xlabel={x},
			ylabel={y},
			enlargelimits={abs=20pt},
			view = {10}{10}
			]


			\addplot3 [
    				domain=-2:2,
    				domain y = -2:2,
    				only marks,
    				mark size=0.5pt,
    				gray] {0};
    			\draw (-2, -2, 0) node[left] {\textcolor{gray}{$\mathcal{M}$}}	;
    			\addplot3 [domain=-2:2, y domain = 0:0, color=gray, samples=5] (x,-2,{0});
    			\addplot3 [domain=-2:2, y domain = 0:0, color=gray, samples=5] (x,2,{0});
    			\addplot3 [domain=-2:2, y domain = 0:0, color=gray, samples=5] (-2,x,{0});
    			\addplot3 [domain=-2:2, y domain = 0:0, color=gray, samples=5] (2,x,{0});
    			
    				\addplot3 [
    				domain=-2:2,
    				domain y = -2:2,
    				only marks,
    				mark size=0.5pt,
    				mesh,
    				draw = black] {0.5 -0.03*( x - y ) + .01*rand};
    			\draw (-2, -2, 0.5) node[left] {$b$}	;
    				
    				\addplot3 [
    				domain=-2:2,
    				domain y = -2:2,
    				surf, 
    				color=blue, 
    				opacity=0.4,
    				faceted color=blue!40] {1};
          	\draw (-2, -2, 1.1) node[left] {\textcolor{blue}{$\eta$}}	;
          		
			\draw[black, thick, ->] (2.1, -1.1, 0) -- (2.1+.5, -1.1, 0) node[right] {$x$} ;
			\draw[black, thick, ->] (2.1, -1.1, 0) -- (2.1, -1.1+1, 0) node[right] {$y$} ;	
          		
		\end{axis}
		\end{tikzpicture}
		\label{fig.water at rest_lake}
	}
	\hfill
	\subfigure[water at rest equilibrium in spherical coordinates]{
		\begin{tikzpicture}[scale=0.8]
			
			\begin{axis}[
			hide axis,
			xlabel={x},
			ylabel={y},
			enlargelimits={abs=20pt},
			view = {10}{10}
			]


			\addplot3 [
    				domain=-1:1,
    				domain y = -1:1,
    				only marks,
    				mark size=0.5pt,
    				gray] {sqrt(4^2 - x^2 - y^2) };
    			\draw (1, -1, 3.8) node[right] {\textcolor{gray}{$\mathcal{M}$}}	;
    			\addplot3 [domain=-1:1, y domain = 0:0, color=gray, samples=20] (x,-1,{sqrt(4^2 - x^2 - 1^2)});
    			\addplot3 [domain=-1:1, y domain = 0:0, color=gray, samples=20] (x,1,{sqrt(4^2 - x^2 - 1^2)});
    			\addplot3 [domain=-1:1, y domain = 0:0, color=gray, samples=20] (-1,x,{sqrt(4^2 - x^2 - 1^2)});
    			\addplot3 [domain=-1:1, y domain = 0:0, color=gray, samples=20] (1,x,{sqrt(4^2 - x^2 - 1^2)});
    			
    				\addplot3 [
    				domain=-1:1,
    				domain y = -1:1,
    				only marks,
    				mark size=0.5pt,
    				mesh,
    				draw = black] {sqrt(4^2 - x^2 - y^2) + .1 + .008*rand};
    			\draw (1, 1, 3.9) node[right] {$b$}	;
    			
    				\addplot3 [
    				domain=-1:1,
    				domain y = -1:1,
    				surf, 
    				color=blue, 
    				opacity=0.4,
    				faceted color=blue!40] {sqrt(4^2 - x^2 - y^2) + .3};
          	\draw (1, 1, 4.1) node[right] {\textcolor{blue}{$\eta$}}	;
    			
    			\draw[->, thick] plot [smooth, domain=0:.35, samples=5] ( {-1+1*\x}, {-1}, {sqrt(4^2 - (-1+1*\x)^2 - 1^2)} );
			\draw[->, thick] plot [smooth, domain=0:.35, samples=5] ( {-1}, {-1+1*\x}, {sqrt(4^2 - (-1+1*\x)^2 - 1^2)} );
			
			\draw ( -.65, -1, 3.818 ) node[below] {$\theta$};
			\draw ( -1, -.65, 3.818 ) node[left] {$\varphi$};
          		
		\end{axis}
		\end{tikzpicture}
		\label{fig.water at rest_ocean}
	}
	\caption{Sketches of water at rest equilibria in Cartesian coordinates (left) and spherical coordinates (right). The meshed surface represents the bottom topography and the blue surface is the free surface. The case in the left fits the equilibrium of the fluid in a lake, with vertical gravitational field. Consequently the free surface is constant for any couple $(x,y)$ of the domain. On the right it is depicted the water at rest solution for an ocean, with the gravitational field pointing towards the center of the planet. In that case, the free surface is constant with respect to the longitudinal and latitudinal angles $(\theta, \varphi)$.}
	\label{fig.water at rest}
\end{figure}

\begin{remark}
	For equilibrium of Prop.~\ref{prop.wr}, condition~\eqref{eq.WRfs} on constant free surface $\eta$ has to be intended with respect to the frame of general covariant coordinates $\x = (x^1, x^2)$. This definition, thus, is still related to the shape of the gravitational field. For example, let us consider two cases of a lake and of an ocean. In the former, the gravitational field can be approximated to be parallel to the vertical direction, in the latter the gravitational field is defined addressing to the center of a spheric planet of radius $R$. Consequently the natural frames of the two cases are Cartesian ($\x = (x,y)$, with $x$ and $y$ coordinates measuring lengths) and spherical ($\x = (\theta, \varphi)$, with $\theta$ and $\varphi$ latitudinal and longitudinal angles), respectively. For the lake, condition~\eqref{eq.WRfs} means that the free surface is constant with respect to Cartesian coordinates $(x,y)$ (see Fig.~\ref{fig.water at rest_lake}), while for the ocean it is intended constant with respect to the angles $(\theta, \varphi)$ but not with respect to the related Cartesian coordinates given by $(R \cos(\theta) \cos(\varphi), R \sin(\theta) \cos(\varphi))$ (see Fig.~\ref{fig.water at rest_ocean}).
\end{remark}

\section{Derivation and final formulation of the novel SW governing model in general covariant coordinates with dependence only on the metric} \label{sec.ourmodel}

In this section, we reformulate the original shallow water model~\eqref{eq.baldauf} in order to make it independent of the physical general covariant reference frame. In particular, our objective consists in giving to the system the metric as input data such that it automatically evolves the physical variables and computes the geometrical features of the manifolds by the only information provided by the metric itself.

\medskip

In order to numerically solve problem~\eqref{eq.baldauf1}, we have re-written it in the first-order hyperbolic form
\begin{equation} \label{eq.hyp}
	\Q_t + \nabla \cdot \F( \Q ) + \B( \Q ) \cdot \nabla \Q = \boldsymbol{0}, \quad \x = (x^1, x^2) \in \Omega \subset \R^d, \quad t \in [0, T],
\end{equation}
where $\x$ is the vector space in the computational domain $\Omega$ of dimension $d = 1,2$ and $t$ is time from 0 to final time $T > 0$. The state $\Q$ (collecting the conserved variables of the equations) is defined in the space of admissible states $\Omega_\Q \in \R^m$. The nonlinear flux tensor and the nonlinear matrix of nonconservative terms are $\F(\Q)$ and $\B(\Q)$, respectively. Their definitions strictly depend on the space dimension $d$. In particular, when $\Omega$ is a two-dimensional set, they read $\F(\Q) = ( \f_1(\Q), \f_2(\Q) )$ and $\B(\Q) = ( \B_1(\Q), \B_2(\Q) )$, with $\f_1, \f_2 \in \R^m$ and $\B_1, \B_2 \in \R^{m \times m}$. Otherwise, they reduce to $\F(\Q) = \f_i(\Q)$ and $\B(\Q) = \B_i(\Q)$ if PDE~\eqref{eq.hyp} is solved along one principal direction led by $x_i$, with $i = 1,2$. System~\eqref{eq.hyp} is hyperbolic if, for any nonzero directions $\n = (n^1, n^2)$, matrix $\mathbf{A} = (\partial \mathbf{F} / \partial \mathbf{Q} + \mathbf{B}(\mathbf{Q}) ) \cdot \mathbf{n} $ has a full real spectrum and its eigenvectors define a complete basis for $\R^m$.

\subsection{Derivation of the governing model} \label{subsec.formulation}
System~\eqref{eq.hyp} is built by splitting all derivatives of the physical quantities from the derivatives of metric terms, as done, for example, in step~\eqref{eq.baldauf1_momWR} of the proof of Prop.~\ref{prop.wr}.\\
\indent In order to consider a full covariant system, indexes of contravariant metric tensor $\gamma^{ij}$ have to be lowered for obtaining its covariant representation $\gamma_{ij}$. The former tensor is the inverse of the latter one and both are symmetric in~$\R^{2 \times 2}$. Let $\alpha : \{ 1, 2 \}^2 \rightarrow \{ 1, 2 \}^2$ be a discrete map defined as $\alpha(i,j) = (3 - i, 3 - j)$. Given an index equal to either 1 or 2, map $\alpha$ exchanges its value with the only other possible one (e.g. $\alpha(1,2) = (2,1)$; $\alpha(2,2) = (1,1)$). Moreover, tensor $G_{ij} = 2 \delta_{ij} -1$ is defined. Since it holds
\begin{equation} \label{eq.metric_con_cov}
	[\gamma^{ij}] = \frac{1}{\gamma} \begin{bmatrix}
		\gamma_{22} & -\gamma_{12} \\
		-\gamma_{12} & \gamma_{11}
	\end{bmatrix},
\end{equation}
with $\gamma = \det[\gamma_{ij}] = \gamma_{11} \gamma_{22} - \gamma_{12}^2$, expression~\eqref{eq.metric_con_cov} in index representation is
\begin{equation} \label{eq.metric_con_cov_ind}
	\gamma^{ij} = \frac{1}{\gamma} G_{(ij)} \gamma_{(\alpha(i,j))}.
\end{equation}
After having applied the chain rule and due to~\eqref{eq.metric_con_cov_ind}, the flux term in the momentum equilibrium~\eqref{eq.baldauf1_mom} is rewritten as
\begin{equation} \label{eq.baldauf_flux}
	\begin{aligned}
		\partialfra{j} T^{ij} &= \partialfra{j} \bigg( \frac{m^i m^j}{h} \bigg) + g h \gamma^{ij} \partialfra{j} h + \frac{1}{2} g h ^2 \partialfra{j} \gamma^{ij} \\
		&= \partialfra{j} \bigg( \frac{m^i m^j}{h} \bigg) + \frac{g h}{\gamma} G_{(ij)} \gamma_{(\alpha(i,j))} \partialfra{j}h + \frac{1}{2} g h^2 \partialfra{j} \bigg( \frac{1}{\gamma} G_{(ij)} \gamma_{(\alpha(i,j))} \bigg).
	\end{aligned}
\end{equation}
For a generic function $f \in C^1(\Omega)$, it holds
\begin{equation} \label{eq.propef}
	\partialfra{k} \frac{f}{\gamma} = \frac{1}{\gamma} \partialfra{k} f - \frac{f}{\gamma^2} \gamma_{22} \partialfra{k} \gamma_{11} + \frac{2 f}{\gamma^2} \gamma_{12} \partialfra{k} \gamma_{12} - \frac{f}{\gamma^2} \gamma_{11} \partialfra{k} \gamma_{22}, \quad k = 1,2.
\end{equation}
Consequently, the last term in~\eqref{eq.baldauf_flux} turns into
\begin{equation} \label{eq.NCPgamma}
	\frac{1}{2} g h^2 \partialfra{j} \bigg( \frac{G_{(ij)}}{\gamma} \gamma_{(\alpha(i,j))} \bigg) = \frac{g h^2}{2 \gamma} \bigg[ \partialfra{j} (G_{(ij)} \gamma_{(\alpha(ij))}) - \gamma_{22} G_{(ij)} \gamma_{(\alpha(ij))} \partialfra{j} \gamma_{11} + 2 \gamma_{12} G_{(ij)} \gamma_{(\alpha(ij))} \partialfra{j} \gamma_{12} - G_{(ij)} \gamma_{(\alpha(ij))} \gamma_{22} \partialfra{j} \gamma_{22} \bigg].
\end{equation}
Relation~\eqref{eq.metric_con_cov_ind} also allows to rewrite the source term in~\eqref{eq.baldauf1_mom} as
\begin{equation} \label{eq.NCPb}
	S^i = g h \gamma^{ij} \partialfra{j}b = \frac{g h}{\gamma} G_{(ij)} \gamma_{(\alpha(ij))} \partialfra{j}b.
\end{equation}
Finally, it remains to treat the terms involving Christoffel symbols in~\eqref{eq.baldauf1}. It is generally valid the contraction property on repeated indexes for Christoffel symbols of second kind:
\begin{equation} \label{eq.contraction}
	\Gamma^j_{jk} = \frac{1}{\sqrt{\gamma}} \partialfra{j} \sqrt{\gamma}.
\end{equation}
Moreover, for the derivative of the square root of the determinant of the metric it holds
\begin{equation} \label{eq.metricder}
	\partialfra{k} \sqrt{\gamma} = \frac{1}{2 \sqrt{\gamma}} \bigg( \gamma_{22} \partialfra{k} \gamma_{11} - 2 \gamma_{12} \partialfra{k} \gamma_{12} + \gamma_{11} \partialfra{k} \gamma_{22} \bigg).
\end{equation} 
In~\eqref{eq.baldauf1}, it follows that terms involving Christoffel symbols with repeated indexes explicitly read
\begin{equation} \label{eq.NCPGamma_rep}
	\begin{aligned}
		\Gamma^j_{jk} m^k &= \frac{m^k}{2 \gamma} \bigg( \gamma_{22} \partialfra{k} \gamma_{11} - 2 \gamma_{12} \partialfra{k} \gamma_{12} + \gamma_{11} \partialfra{k} \gamma_{22} \bigg), \\
		\Gamma^j_{jk} T^{ik} &= \frac{T^{ik.}}{2 \gamma} \bigg( \gamma_{22} \partialfra{k} \gamma_{11} - 2 \gamma_{12} \partialfra{k} \gamma_{12} + \gamma_{11} \partialfra{k} \gamma_{22} \bigg),
	\end{aligned}
\end{equation}
in the mass and in the momentum balances, respectively. Concerning the terms involving the Christoffel symbols with not repeated indexes, by using property~\eqref{eq.propef} with $f \equiv 1$ to be applied in definition~\eqref{eq.christoffel}, one can write
\begin{subequations} \label{eq.NCPGamma}	
	\begin{equation} \label{eq.NCPGamma1}
		\begin{aligned}
			\Gamma^1_{jk} T^{kj} = &\frac{T^{11}}{\gamma} \gamma_{22} \partialfra{1}\gamma_{11} - 2 \frac{T^{11}}{\gamma} \gamma_{12} \partialfra{1} \gamma_{12} + \frac{1}{2 \gamma} ( T^{11} \gamma_{11} - 2 T^{12} \gamma_{12} - D^{22} \gamma_{22}) \partialfra{1} \gamma_{22} \\
			&+ \frac{1}{2 \gamma} ( T^{11} \gamma_{12} + 3 T^{12} \gamma_{22}) \partialfra{2} \gamma_{11} + \frac{1}{\gamma} (T^{22} \gamma_{22} - T^{12} \gamma_{12}) \partialfra{2} \gamma_{12} + \frac{1}{2\gamma}(T^{12} \gamma_{11} - T^{22} \gamma_{12}) \partialfra{2} \gamma_{22},
		\end{aligned}
	\end{equation}
	\text{and}
	\begin{equation} \label{eq.NCPGamma2}
		\begin{aligned}
			\Gamma^2_{jk} T^{kj} = &\frac{1}{2\gamma} (T^{12} \gamma_{22} - T^{11} \gamma_{12}) \partialfra{1} \gamma_{11} + \frac{1}{\gamma} (T^{11} \gamma_{11} - T^{12} \gamma_{12}) \partialfra{1} \gamma_{12} + \frac{1}{2\gamma} (T^{22} \gamma_{12} + 3 T^{12} \gamma_{11} \partialfra{1} \gamma_{22}) \\
			&+ \frac{1}{2\gamma} ( T^{22} \gamma_{22} - T^{11} \gamma_{11} - 2 T^{12} \gamma_{12}) \partialfra{2} \gamma_{11} - 2 \frac{T^{22}}{\gamma} \gamma_{12} \partialfra{2} \gamma_{12} + \frac{T^{22}}{\gamma} \gamma_{11} \partialfra{2} \gamma_{11}.
		\end{aligned}
	\end{equation}
\end{subequations}


\subsection{Final formulation of the governing model}
It is now possible to define all the components in~\eqref{eq.hyp}. The state $\Q$ is
\begin{equation} \label{eq.Q}
	\Q = [h, m^1, m^2, b, \gamma_{11}, \gamma_{12}, \gamma_{22}]^T,
\end{equation} 
where the first four entries are the physical variables and the remaining three entries completely describe the covariant metric tensor (that is symmetric and, thus, completely defined by the diagonal terms $\gamma_{11}$ and $\gamma_{22}$ and the extra-diagonal term $\gamma_{12} = \gamma_{21}$). Consequently, the dimension $m$ of the space $\Omega_\Q$ of the admissible states is 7, with the first three variables to be evolved, the bathymetry in fourth position, and the last three components being the metric quantities in covariant representation. Flux $\F(\Q)$ is purely kinetic and it is defined by\footnote{Expressions $\mathbf{0}_a$ and $\mathbf{O}_{b \times c}$, with $a, b, c \in \mathbb{N}$, mean the null vector in $\R^a$ and the null matrix in $\R^{b \times c}$, respectively.}
\begin{equation} \label{eq.F}
	\f_1 = \begin{bmatrix}
		m^1 \\[3pt]
		\frac{(m^1)^2}{h} \\[3pt]
		\frac{m^1 m^2}{h} \\[3pt]
		\mathbf{0}_4
	\end{bmatrix}, \quad \text{and} \quad \f_2 = \begin{bmatrix}
		m^2 \\[3pt]
		\frac{m^1 m^2}{h} \\[3pt]
		\frac{(m^2)^2}{h} \\[3pt]
		\mathbf{0}_4
	\end{bmatrix}.
\end{equation}
Finally, matrices collecting the nonconservative parts can be split into two blocks
\begin{equation} \label{eq.B}
	\B_i(\Q) = \B_i^\eta(\Q) + \B_i^{met}(\Q), \quad i = 1,2, 
\end{equation} 
where $\B_i^\eta(\Q)$ is the block obtained from the derivatives of $h$ in~\eqref{eq.baldauf_flux} and $b$ in~\eqref{eq.NCPb} and $\B_i^{met}(\Q)$ is the block gathering the derivatives of the covariant metric components $\gamma_{ij}$ in~\eqref{eq.NCPgamma},~\eqref{eq.NCPGamma_rep} and~\eqref{eq.NCPGamma}. They read
\begin{subequations} \label{eq.Bexpr}
	\begin{equation} \label{eq.Beta}
		\B_1^\eta = \frac{1}{\gamma} \left[\begin{array}{c|c}
			\begin{matrix}
				0 & 0 & 0 & 0 \\[3pt]
				g h \gamma_{22} & 0 & 0 & g h \gamma_{22} \\[3pt]
				-g h \gamma_{12} & 0 & 0 & -g h \gamma_{12} \\[3pt]
				0 & 0 & 0 & 0 
			\end{matrix} & \mathbf{O}_{4 \times 3} \\[3pt] \hline
			\mathbf{O}_{3 \times 4} & \mathbf{O}_{3 \times 3}
		\end{array} \right], \quad \B_2^\eta = \frac{1}{\gamma} \left[ \begin{array}{c|c}
			\begin{matrix}
				0 & 0 & 0 & 0 \\[3pt]
				-g h \gamma_{12} & 0 & 0 & -g h \gamma_{12} \\[3pt]
				g h \gamma_{11} & 0 & 0 & g h \gamma_{11} \\[3pt]
				0 & 0 & 0 & 0 
			\end{matrix}	& \mathbf{O}_{4 \times 3} \\[3pt] \hline
			\mathbf{O}_{3 \times 4} & \mathbf{O}_{3 \times 3}		
		\end{array}	\right],
	\end{equation}
	\text{and}
	\begin{equation} \label{eq.Bmet}
	\begin{aligned}
		\B_1^{met} &= \frac{1}{\gamma} \left[ \begin{array}{c|c}
			\mathbf{O}_{3 \times 4} & \begin{matrix}
				\frac{1}{2} m^1 \gamma_{22} & -m^1 \gamma_{12} & \frac{1}{2} m^1 \gamma_{11} \\[3pt]
				\frac{(m^1)^2}{h} \gamma_{22} & -2 \frac{(m^1)^2}{h} \gamma_{12} & -\frac{1}{2 h} [-(m^1)^2 \gamma_{11} + 2 m^1 m^2 \gamma_{12} + (m^2)^2 \gamma_{22}] \\[3pt] 
				\frac{m^1}{2 h} (-m^1 \gamma_{12} + m^2 \gamma_{22}) & \frac{m^1}{h} (m^1 \gamma_{11} - m^2 \gamma_{12}) & \frac{m^2}{2 h} (3 m^1 \gamma_{11} + m^2 \gamma_{12})
			\end{matrix} \\[3pt] \hline
			\mathbf{O}_{4 \times 4} & \mathbf{O}_{4 \times 3}
		\end{array} \right], \\
		\B_2^{met} &= \frac{1}{\gamma} \left[ \begin{array}{c|c}
			\mathbf{O}_{3 \times 4} & \begin{matrix}
				\frac{1}{2} m^2 \gamma_{22} & -m^2 \gamma_{12} & \frac{1}{2} m^2 \gamma_{11} \\[3pt]
				\frac{m^1}{2 h} (3 m^2 \gamma_{11} + m^1 \gamma_{12}) & \frac{m^2}{h} (m^2 \gamma_{11} - m^1 \gamma_{12}) & \frac{m^2}{2 h} (-m^2 \gamma_{12} + m^1 \gamma_{22}) \\[3pt] 
				-\frac{1}{2 h} [-(m^2)^2 \gamma_{11} + 2 m^1 m^2 \gamma_{12} + (m^1)^2 \gamma_{22}] & -2 \frac{(m^2)^2}{h} \gamma_{12} & \frac{(m^2)^2}{h} \gamma_{22}
			\end{matrix} \\[3pt] \hline
			\mathbf{O}_{4 \times 4} & \mathbf{O}_{4 \times 3}
		\end{array} \right],
	\end{aligned}
	\end{equation}
\end{subequations} 
respectively.

\section{Numerical method} \label{sec.method}

Once system~\eqref{eq.baldauf1} is re-written in hyperbolic form~\eqref{eq.hyp} as previously explained, it is possible to evolve the solution at discrete level by employing a second-order MUSCL-Hancock scheme. For now it is assumed that the preserved equilibrium is smooth, i.e. it does not present any discontinuity at the element interfaces. However, for the construction of the WB scheme for water at rest equilibrium~\eqref{eq.WR} it will be possible to relax this hypothesis.

\medskip

The numerical method is presented in two dimensions with the computational domain $\Omega$ discretized by a mesh of Voronoi-type cells. The scheme is formally the same when the dimension $d$ is either one or two. In particular, when the test case is one-dimensional, instead of having a polygonal tassellation, the cells reduce to intervals of equal length. A consequent simplification of the scheme is resumed in Sec.~\ref{subsec.1D}. Concerning the two-dimensional case, the advantages in employing a Voronoi-type discretization instead of more classical triangulations are largely analyzed in~\cite{boscheriAFE2022}. 

Let $\mathcal{T}_N$, with $N \in \N$, be a proper partition of the computational domain $\Omega$. Any cell $\Omega_k \in \mathcal{T}_N$, for $k  = 1, \dots, N$, is built such that $\Omega_k \cap \Omega_l = \emptyset$, for any $k \neq l$, and $\bigcup_{k = 1}^N \Omega_k \equiv \Omega$. The stencil set of neighbors sharing one edge with cell $\Omega_k$ is denoted by $\mathcal{S}_k$. Consequently, let $\E_k$ be the set of edges of cell $\Omega_k$; in particular edge $e_{kl} \in \E_k$ is the edge shared by cell $\Omega_k$ and its neighboring cell $\Omega_l \in \mathcal{S}_k$ (i.e. $\overline{\Omega_k \cap \Omega_l} = e_{kl}$). Let $t^n$ be the $n$-th time instance such that $t^{n+1} = t^n + \Delta t$, with $\Delta t$ the time step. Finally, the cell average of the state over cell $\Omega_k$ at time $t^n$ is denoted by
\begin{equation} \label{eq.cellave}
	\Q_k^n = \frac{1}{|\Omega_k|} \int_{\Omega_k} \Q(\x, t^n) \, \diff{\Omega}.
\end{equation}

\subsection{Second order finite volume MUSCL-Hancock-type scheme} \label{subsec.MH}

The standard MUSCL-Hancock-type scheme, originally introduced in~\cite{van1974towards} and pedagogically presented in~\cite{toro2013riemann}, involves a local polynomial reconstruction $\q_k(\x, t)$ in space and time, which is then used to numerically approximate the fluxes and the jumps at the interfaces in the evolutionary scheme at finite volumes. This reconstruction is defined by
\begin{equation} \label{eq.local_reconstruction}
	\q_k(\x, t) = \mathbf{w}_k(\x) + \partial_t \Q_k (t - t^n) = \Q_k^n + \nabla \Q_k (\x - \x_k) + \partial_t \Q_k (t - t^n), \quad (\x, t) \in \Omega_k \times (t^n, t^{n+1}),
\end{equation}
where $\x_k$ is the center of mass of cell $\Omega_k$ and vectors $\nabla \Q_k$ and $\partial_t \Q_k$ in $\R^m$ are the local vector polynomial coefficients to be found.\\
Concerning the spatial reconstruction $\mathbf{w}_k$, it comes from imposing the integral conservation to be preserved in any cell of the stencil, i.e.
\begin{equation} \label{eq.integral_conservation}
	\frac{1}{|\Omega_l|} \int_{\Omega_l} \mathbf{w}(\x) \, \diff{\Omega} = \Q_l^n, \quad \forall \Omega_l \in \mathcal{S}_k;
\end{equation}
nevertheless, since condition~\eqref{eq.integral_conservation} is usually over-determined, it is strongly imposed that the constraint is exactly valid only for cell $\Omega_k$, consequently, equation~\eqref{eq.integral_conservation} is rewritten as
\begin{equation} \label{eq.DeltaQk}
	\frac{1}{|\Omega_l|} \int_{\Omega_l} \nabla \Q_k (\x - \x_k) \, \diff{\Omega} = \Q_l^n - \Q_k^n, \quad \forall \Omega_l \in \mathcal{S}_k.
\end{equation}
Then, Problem~\eqref{eq.DeltaQk} is solved in the sense of least-squares. This provides a non-limited slope $\nabla \Q_k$. \RIIcolor{ In order not to create new extrema in the spatial reconstruction process, a slope limiter $\Phi_k |_{\nabla \Q_k}$ is employed~\cite{toro2013riemann, van1974towards, van1979towards}. In particular, let $w_k^{\mu, \max}$ and $w_k^{\mu, \min}$ be the maximum and minimum value of the $\mu$-th entry of a vector $\mathbf{w}_k \in \R^m$, respectively, with $\mu = 1, \dots, m$. For any vertex $\mathbf{v}_p^{kl}$, with $p = 1,2$, of edge $e_{kl} \in \mathcal{E}_k$, a local element-wise and component-wise slope-limiter $\Phi_{kl,p}^\mu |_{\mathbf{w}}$~\cite{barth1989design} applied to vector $\mathbf{w}$ is defined as
\begin{equation} \label{eq.local_slopelimiter}
	\Phi_{kl,p}^\mu \bigg|_{\mathbf{w}} = \left\{ \begin{array}{ll}
		\min \left(1, \frac{w_k^{\mu, \max} - w_k^{\mu}}{w_k^\mu(\mathbf{v}_p^{kl}) - w_k^{\mu}  } \right), & \text{if } w_k^\mu(\mathbf{v}_p^{kl}) - w_k^{\mu}  > 0, \\
		\min \left(1, \frac{w_k^{\mu, \min} - w_k^{\mu}}{w_k^\mu(\mathbf{v}_p^{kl}) - w_k^{\mu}  } \right), & \text{if } w_k^\mu(\mathbf{v}_p^{kl}) - w_k^{\mu}  < 0. \\
		1, & \text{if } w_k^\mu(\mathbf{v}_p^{kl}) - w_k^{\mu}  = 0.
	\end{array} \right.
\end{equation} 
Thus $\Phi_k |_{\mathbf{w}} = [ \min_{e_{kl}  \in \mathcal{E}_k} ( \min_{\mathbf{v}_p \in e_{kl}} ( \Phi_{kl,p}^\mu |_{\mathbf{w}} )) ]_{\mu = 1}^m \in \R^m$. By abuse of notation, we still write $\nabla \Q_k$ for the limited slope of the local spatial reconstruction instead of $\Phi_k |_{\nabla \Q_k}$}. Regarding the time coefficient $\partial_t \Q_k$ in~\eqref{eq.local_reconstruction}, it is directly computed by discretely integrating the hyperbolic formulation~\eqref{eq.hyp} along the boundary $\partial \Omega_k = \bigcup_{e_{kl} \in \mathcal{E}_k} e_{kl}$ of cell $\Omega_k$ as follows
\begin{equation} \label{eq.Qt}
	\partial_t \Q_k = - \sum_{e_{kl} \in \mathcal{E}_k } \bigg( | e_{kl} | \F(\q_k(\x_{e_{kl}}, t^n) ) \cdot \n_{e_{kl}} \bigg) - \B(\Q_k^n) \cdot \nabla \Q_k,  
\end{equation}
where $\x_{e_{kl}}$ is the midpoint of edge $e_{kl}$ whose outward unit normal vector is $\n_{e_{kl}}$.\\
\indent When the local spatial reconstruction~\eqref{eq.local_reconstruction} is computed, it is possible to use it in the cell average evolution from time $t^n$ to $t^{n+1}$. In particular, for standard MUSCL-Hancock finite volume strategy, the scheme reads
\begin{equation} \label{eq.FV}
		\Q_k^{n+1} = \Q_k^n - \frac{\Delta t}{| \Omega_k |} \sum_{e_{kl} \in \mathcal{E}_k } |e_{kl}| \boldsymbol{\mathcal{F}}_{e_{kl}}( \q_{e_{kl}}^-, \q_{e_{kl}}^+) \cdot \n_{e_{kl}} - \frac{\Delta t}{| \Omega_k |} \sum_{e_{kl} \in \mathcal{E}_k } \boldsymbol{\mathcal{D}}_{e_{kl}}( \q_{e_{kl}}^-, \q_{e_{kl}}^+) \RIIcolor{\cdot \n_{e_{kl}}} - \Delta t \, \B(\q_k^{n + \frac{1}{2}}) \cdot \nabla \Q_k,
\end{equation}
where $\boldsymbol{\mathcal{F}}_{e_{kl}}$ and $\boldsymbol{\mathcal{D}}_{e_{kl}}$ represent the numerical flux and the nonconservative jump at the interface $e_{kl}$, respectively. A Rusanov approach is introduced in order to approximate the flux, i.e.
\begin{equation} \label{eq.rusanov}
	\boldsymbol{\mathcal{F}}_e (\q^-, \q^+) \cdot \n_e = \frac{1}{2} \left( \F(\q^+) + \F(\q^-) \right) \cdot \n_e - \frac{1}{2} s_{\max} \tilde{\mathbf{I}} (\q^- - \q^+),
\end{equation}
where $s_{\max}$ is the maximum eigenvalue in the union of spectra of both matrices $\mathbf{A}(\q^+)$ and $\mathbf{A}(\q^-)$. In~\eqref{eq.rusanov}, matrix $\tilde{\mathbf{I}} \in \R^{m \times m}$ is a modified identity matrix whose diagonal reads $\text{diag}(\tilde{\mathbf{I}}) = [1, 1, 1, 0, 0, 0, 0]$. In particular, when it is known \textit{a priori} that a variable is not evolving in time (as it happens for the bathymetry and all the metric quantities), its flux does not need to be stabilized; consequently it is simply approximated as the mean value of the fluxes evaluated from outside and inside the interface. \RIIcolor{ For the jumps at the interface, a Lipshitz-continuous path $\boldsymbol{\Psi}(\q^-, \q^+; \tau)$, with $\tau \in [0, 1]$, is employed~\cite{dal1995definition} for defining
\begin{equation} \label{eq.jump}
	\boldsymbol{\mathcal{D}}_e (\q^-, \q^+) = \frac{1}{2} \int_0^1  \B( \boldsymbol{\Psi}(\q^-, \q^+; \tau) )  \, \diff{\tau} \, \cdot (\q^+ - \q^-).
\end{equation} 
The conservative path $\boldsymbol{\Psi}(\q^-, \q^+; \tau)$ has to be built such that it assumes values $\q^-$ and $\q^+$ when $\tau = 0$ and $\tau = 1$, respectively. In order to ensure a second-order convergence, it is sufficient to have a simple straight line segment path 
\begin{equation} \label{eq.Psi}
	\boldsymbol{\Psi}(\q^-, \q^+; \tau) = \q^- + \tau (\q^+ - \q^-).
\end{equation} 
In the simulations of this paper, the path $\Psi$ is independent of normal direction $\n_e$ of the edge. It is possible to define a path depending on the direction and the interested reader is referred to~\cite{castro2010some}.}
In the scheme~\eqref{eq.FV} the following notations are used:
\begin{equation} \label{eq.FVnotation}
	\q_{e_{kl}}^- = \q_k(\x_{e_{kl}}, t^{n + \frac{1}{2}}), \quad \q_{e_{kl}}^- = \q_l(\x_{e_{kl}}, t^{n + \frac{1}{2}}), \quad \q_k^{n + \frac{1}{2}} = \q_k(\x_k, t^{n + \frac{1}{2}}),
\end{equation}
with $t^{n + \frac{1}{2}} = t^n + \frac{1}{2} \Delta t$ and $\q$ defined in~\eqref{eq.local_reconstruction}.\\

\subsection{MUSCL-Hancock scheme in one dimension} \label{subsec.1D}
In one dimension, the MUSCL-Hancock scheme previously presented largely simplifies. Let $\xi$ denote the spatial direction $x^1$ or $x^2$ along which the problem is integrated. The domain $\Omega = [\xi_L, \xi_R]$ is considered to be covered by $N$ intervals $\Omega_k = [\xi_{k - \frac{1}{2}}, \xi_{k + \frac{1}{2}}]$ of length $\Delta \xi = (\xi_R - \xi_L)/N$.\\
The spatial reconstruction $\mathbf{w}_k$ in the space-time polynomial~\eqref{eq.local_reconstruction} reads
\begin{equation} \label{eq.spatial_local_reconstruction1D}
	\mathbf{w}_k(\xi) = \Q_k^n + \partial_\xi \Q_k (\xi - \xi_k), \quad \xi \in \Omega_k,
\end{equation}
where $\xi_k$ is the middle point of cell $\Omega_k$ and the limited spatial derivative $\partial_\xi \Q_k$ is computed by employing the $\minmod$ function between the first order approximation of the derivatives from the right and from the left of the cell. Namely,
\begin{equation} \label{eq.Qxi}
	\partial_\xi \Q_k = \minmod\left( \frac{\Delta \Q_{k + \frac{1}{2}}}{\Delta \xi}, \frac{\Delta \Q_{k - \frac{1}{2}}}{\Delta \xi} \right),
\end{equation}
with differences $\Delta \Q_{k + \frac{1}{2}} = \Q_{k+1}^n - \Q_k^n$ and $\Delta \Q_{k - \frac{1}{2}} = \Q_k^n - \Q_{k -1}^n$ at the right and at the left of the cell, respectively, and the $\minmod$ function of two reals $a$ and $b$ defined as
\begin{equation} \label{eq.minmod}
	\minmod(a, b) = \left\{ \begin{array}{ll}
		0, & \text{if } ab \leq 0, \\
		a, & \text{if } |a| < |b|, \\
		b, & \text{if } |a| \geq |b|.
	\end{array}	 \right. 
\end{equation}
The time term $\partial_t \Q_k$ in~\eqref{eq.local_reconstruction} is directly given by the hyperbolic formulation~\eqref{eq.hyp} as
\begin{equation} \label{eq.Qt1D}
	\partial_t \Q_k = - \frac{ \F(\Q_k^n + \frac{1}{2} \Delta\xi \partial_\xi \Q_i^n) - \F(\Q_k^n - \frac{1}{2} \Delta\xi \partial_\xi \Q_i^n) }{ \Delta \xi } - \B(\Q_k^n) \cdot \partial_\xi \Q_k^n,
\end{equation}
where the derivative of the flux term is approximated at the first order from the right to the left extreme of interval $\Omega_k$.\\
Consequently, the finite volume scheme~\eqref{eq.FV} reduces to
\begin{equation} \label{eq.FV1D}
	\begin{aligned}
		\Q_k^{n+1} = \Q_k^n &- \frac{\Delta t}{\Delta \xi} \left( \boldsymbol{\mathcal{F}}_{k + \frac{1}{2}}( \q_{k + \frac{1}{2}}^-, \q_{k + \frac{1}{2}}^+) - \boldsymbol{\mathcal{F}}_{k - \frac{1}{2}}( \q_{k - \frac{1}{2}}^-, \q_{k - \frac{1}{2}}^+) \right)\\ 
		&- \frac{\Delta t}{\Delta \xi} \left( \boldsymbol{\mathcal{D}}_{k + \frac{1}{2}}( \q_{k + \frac{1}{2}}^-, \q_{k + \frac{1}{2}}^+) + \boldsymbol{\mathcal{D}}_{k - \frac{1}{2}}( \q_{k - \frac{1}{2}}^-, \q_{k - \frac{1}{2}}^+) \right)\\ 
		&- \Delta t \, \B(\q_k^{n + \frac{1}{2}}) \cdot \partial_\xi \Q_k.
	\end{aligned}
\end{equation}

\subsection{Time discretization} \label{subsec.time}
Finally, a standard CFL condition is applied in order to determine the time step $\Delta t$. For an unstructured polygonal 2D tessellation the CFL constraint reads as follows~\cite{gaburro2020high}
\begin{equation} \label{eq.CFL}
	\Delta t = \CFL \min_{ \Omega_k \in \mathcal{T}_N } \frac{ | \Omega_k | }{ | \lambda_{k, \max} | \sum_{e_{kl} \in \mathcal{E}_k } | e_{kl} | },
\end{equation}
where $\CFL < 1/d$~\cite{dumbser2008unified} is the Courant-Friederichs-Levy number and $| \lambda_{k, \max} |$ is the absolute value of the maximum eigenvalue of matrix $\mathbf{A}$ over cell $\Omega_k$, and for the 1D case it simplifies to
\begin{equation}
	\Delta t = \textnormal{CFL} \, \min \limits_{\Omega_k \in \mathcal{T}_N } \frac{\Delta \xi}{|\lambda_{k, \max}|}. 
	\label{eq:time step}
\end{equation}

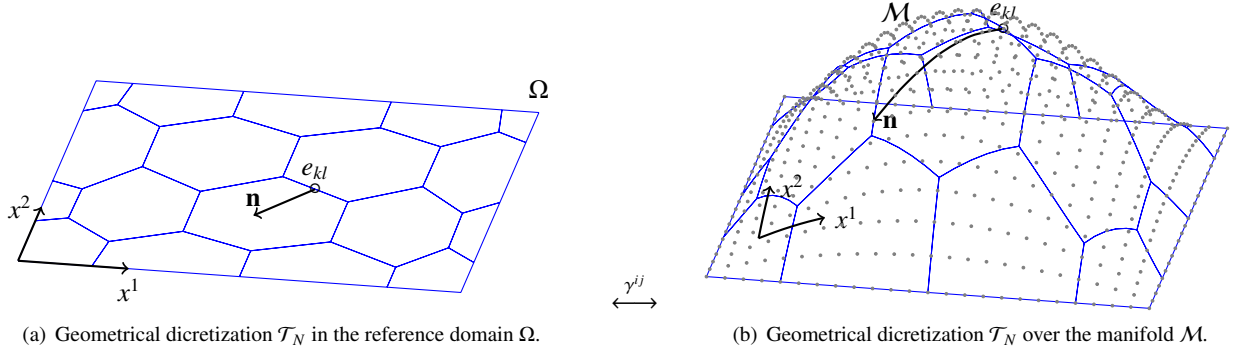
\begin{figure}
	\centering
	\subfigure[Geometrical dicretization $\mathcal{T}_N$ in the reference domain $\Omega$.]{
		\begin{tikzpicture}
		\begin{axis}[
	hide axis,
	enlargelimits={abs=10pt},
	xlabel={x},
	ylabel={y},
    view = {10}{40}
]
	
	\addplot3 [mark=none, blue] table {./data/Test_B_EUL_CFL100-090_Lim23_RefDG2_Edges.dat};
	
	\draw (0.165, 0.0, 0) node[] {$\circ$};
	\draw [->, thick] (0.165, 0.0, 0) node[above] {$e_{kl}$} -- (-0.0550, -0.3300, 0) node[above] {$\n$};

\draw (1, 1, .05) node[above] {$\Omega$};
\draw [->, thick] (-1, -1, 0) -- (-.5, -1, 0) node[below] {$x^1$};
\draw [->, thick] (-1, -1, 0) -- (-1, -.4, 0) node[left] {$x^2$};

\end{axis}
	\end{tikzpicture}
	\label{fig.geometry_ref}
	}
	\hfill
	$ \overset{\gamma^{ij}}{\longleftrightarrow}$ 
	\hfill
	\subfigure[Geometrical dicretization $\mathcal{T}_N$ over the manifold $\mathcal{M}$.] {
		\begin{tikzpicture}
\begin{axis}[
	hide axis,
	xlabel={x},
	ylabel={y},
    view = {10}{40}
]
	\addplot3 [mark=none, blue] table {./data/voronoi_b.dat};

			\addplot3 [
    				domain=-1:1,
    				domain y = -1:1,
    				only marks,
    				mark size=0.5pt,
    				gray] {.05*(x-1)*(x+1)*(y-1)*(y+1) };

\draw (0.165, 0.0, 0.0486) node[] {$\circ$};
\draw[<-, thick] plot [smooth, domain=-.5:0, samples=5] ( {0.165+1*\x}, {1*\x}, {.05*(0.165+1*\x-1)*(0.165+1*\x+1)*(\x-1)*(\x+1)} );
\draw (0.165, 0.0, 0.0486) node[above] {$e_{kl}$};
\draw (-0.3350, -.5, 0.0333) node[right] {$\n$};


\draw[->, thick] plot [smooth, domain=0:.3, samples=5] ( {-.8+1*\x}, {-.8}, {.05*(-.8+1*\x-1)*(-.8+1*\x+1)*(-.8-1)*(-.8+1)} );
\draw[->, thick] plot [smooth, domain=0:.3, samples=5] ( {-.8}, {-.8+1*\x}, {.05*(-.8+1*\x-1)*(-.8+1*\x+1)*(-.8-1)*(-.8+1)} );
\draw (-0.5000, -.8, 0.0135) node[right] {$x^1$};
\draw ( -.8,-0.5000, 0.0135) node[right] {$x^2$};

\draw (-.5, 1, 2e-2) node[above] {$\mathcal{M}$};

\end{axis}
\end{tikzpicture}
\label{fig.geometry_phys}
	}
	\caption{On the left, $\n$ is the unit normal vector to edge $e_{kl}$ in the computational domain $\Omega$. When vector $\n$ is represented in the physical frame over the manifold $\mathcal{M}$ (figure on the right), its length is changed due to the curvature of the manifold. The metric $\gamma^{ij}$ defines the physical length of this vector $\n$ in the computation of time step~\eqref{eq.CFL}.}
	\label{fig.geometry}
\end{figure}

In this context, it is worth analyzing the spectrum $\rho(\mathbf{A})$ of $\mathbf{A}$. For a given direction $\n$, the nonzero components of $\rho(\mathbf{A})$ are
\begin{equation} \label{eq.eigenvalues}
	\{ \mathbf{u} \cdot \n, \mathbf{u} \cdot \n \pm c \},
\end{equation}
with 
\begin{equation} \label{eq.c}
	c = \sqrt{g h} \, \| \n \|_{\gamma^{ij}}, 
\end{equation}
where $\sqrt{g h}$ is the SW wave speed in Cartesian coordinates and $\| \n \|_{\gamma^{ij}} = \sqrt{ \n^T [\gamma^{ij}] \n }$ is the norm of the direction $\n$ induced by the tensor metric $\gamma^{ij}$. This implies that the metric is influencing the spectrum $\rho(\mathbf{A})$ and, consequently, the definition of the time step~\eqref{eq.CFL}. Moreover, the metric is affecting the measure of the length of the direction $\n$, namely its norm. This can be explained by the fact that the spectrum $\rho(\mathbf{A})$ is computed from system~\eqref{eq.hyp} which is already completely defined in a covariant reference system. Thus, although in the reference domain the direction $\n$ has unit length (see Fig.~\ref{fig.geometry_ref}), its transposition on the manifold $\mathcal{M}$ is affected by deformation due to the curvature of the manifold itself (see Fig.~\ref{fig.geometry_phys}). To account for its deformation along the curvature, the spectrum then involves the metric (used to describe the manifold) for the adjusted measurement of the direction $\n$ and finally to calculate the propagation velocity~\eqref{eq.c} of a wave on a manifold for problem~\eqref{eq.hyp}.

\begin{remark} \label{rem.cartesian1}
	When the metric is Cartesian, i.e. $\gamma^{ij} = \delta^{ij}$, velocity~\eqref{eq.c} reduces to the classical value $c = \sqrt{g h}$. As a matter of fact, the curvature of the associated manifold is 0. Thus, the measurements of the direction along the manifold and in the reference system coincide.\\
\end{remark}

\begin{remark} \label{rem.cartesian2}
	From the previous remark, it follows that, in the Cartesian metric, system~\eqref{eq.hyp} has the same eigenvalues of the classical SW. For this reason, given equal initial conditions, boundary conditions and domain discretization, two numerical solutions for classical SW and SW in Cartesian GCC via MUSCL-Hancock scheme have to coincide all over the domain and at any time. 
\end{remark}

\begin{remark} \label{rem.cartesianWR}
	When the solution for hyperbolic system~\eqref{eq.hyp} is steady, time step $\Delta t$ is fixed because no physical variable is evolving in time. In particular, when the steady solution is specialized to be the water at rest solution (i.e. zero velocity), time step in~\eqref{eq.CFL} only depends on the inverse of the characteristic velocity $c$ if the mesh and the $\CFL$ constant are fixed. The local wave speed solves the optimization problem
	\begin{equation} \label{eq.optimization}
		c = \sqrt{g h} \max_{|\n| = 1} \| \n \|_{\gamma^{ij}}
	\end{equation}
	over any possible unit direction $\n \in \R^2$, which is equivalent to solve
	\begin{equation} \label{eq.maxn}
		\max_{|\n| = 1} \n^T [\gamma^{ij}] \n.
	\end{equation} 
	Due to the min-max theorem by Courant-Fischer-Weyl~\cite{teschl2009mathematical, lieb2001analysis}, problem~\eqref{eq.maxn} is solved by the maximum eigenvalue of the metric tensor $\gamma^{ij}$. After some computation, in covariant representation, the solution of problem~\eqref{eq.maxn} explicitly reads
	\begin{equation} \label{eq.maxn_sol}
		\frac{ \gamma_{11} + \gamma_{22} + \sqrt{ (\gamma_{11} - \gamma_{22})^2 + 4 \gamma_{12}^2 } }{2 \gamma}.
	\end{equation}
	Consequently, the maximum eigenvalue all over the domain is
	\begin{equation} \label{eq.maxeigenvalue}
		\lambda_{\max} = \max_{\x \in \Omega} \sqrt{g h \frac{ \gamma_{11} + \gamma_{22} + \sqrt{ (\gamma_{11} - \gamma_{22})^2 + 4 \gamma_{12}^2 } }{2 \gamma}}.
	\end{equation}
	The water at rest equilibrium offers the possibility to analytically quantify the influence of the metric in the definition of time step $\Delta t$. Thus, once the mesh and the $\CFL$ coefficient are set, it is possible to directly compute the time step~\eqref{eq.CFL} through~\eqref{eq.maxeigenvalue}. Given the same initial physical quantities, it is therefore clear that the metric is the only quantity in the system that defines the time step $\Delta t$. Consequently, to reach a given physical time, the number of iterations required varies with respect to the metric (i.e., with respect to the manifold over which the system is solved).
\end{remark}

\section{Well balanced scheme} \label{sec.WB}
This section is split in two parts. First the well-balancing of the MUSCL-Hancock scheme is explained for a generic equilibrium $\Q^E$ of system~\eqref{eq.hyp}. Afterwards, a minimally invasive and simple well balanced approach is introduced for the water at rest equilibrium~\eqref{eq.WR} denoted in this section with $\Q^{WR}$.

\subsection{Well balanced scheme for a general equilibrium} \label{subsec.WB}
According to~\cite{pares2006numerical, castro2006high, castro2008well}, scheme~\eqref{eq.FV} can be well balanced by the knowledge of a smooth~\cite{berberich2021high, grosheintz2019high, castro2020well} equilibrium $\Q^E$. The objective is to preserve the equilibrium.\\
Let $\Q_{e_{kl}}^E = \Q^E(\x_{kl})$ and $\Q_k^E = \Q^E(\x_k)$ be the equilibrium evaluated over the center of mass $\x_k$ of cell $\Omega_k$ and the midpoint $\x_{kl}$ of edge $e_{kl}$. To evaluate equation~\eqref{eq.hyp} at the equilibrium leads to writing the following exact relation (valid by construction)
\begin{equation} \label{eq.PDEatEq}
	\frac{1}{| \Omega_k |} \sum_{e_{kl} \in \mathcal{E}_k} |e_{kl}| \F(\Q_{e_{kl}}^E) \cdot \n_{e_{kl}} + \frac{1}{| \Omega_k |} \int_{\Omega_k} \B(\Q^E) \cdot \nabla \Q^E \, \diff{\Omega} = 0.
\end{equation}
Consequently the well balanced scheme~\eqref{eq.FV} is built by subtracting relation~\eqref{eq.PDEatEq} from the standard finite volume scheme~\eqref{eq.FV}
\begin{equation} \label{eq.FVWB}
		\Q_k^{n+1} = \Q_k^n - \frac{\Delta t}{| \Omega_k |} \sum_{e_{kl} \in \mathcal{E}_k } |e_{kl}| \left( \boldsymbol{\mathcal{F}}_{e_{kl}}( \q_{e_{kl}}^-, \q_{e_{kl}}^+) - \F(\Q_{e_{kl}}^E) \right) \cdot \n_{e_{kl}} - \frac{\Delta t}{| \Omega_k |} \sum_{e_{kl} \in \mathcal{E}_k } \boldsymbol{\mathcal{D}}_{e_{kl}}( \q_{e_{kl}}^-, \q_{e_{kl}}^+) - \Delta t \left( \B(\q_k^{n + \frac{1}{2}}) \cdot \nabla \Q_k - \B(\Q_k^E) \cdot \nabla \Q_k^E \right) .
	\end{equation}
The discrete spatial variation term $\nabla \Q_k$ depends on the fluctuation $\Q_k^f$ of the current state $\Q_k$ with respect to the chosen equilibrium $\Q^E$ and it will be discussed later.\\
In~\cite{gaburro2018well} the authors propose to split the path linking two states $\q^+$ and $\q^-$ in the sum of the path of two equilibria $\q_E^+$ and $\q_E^-$ and the path of two fluctuations $\q_f^+$ and $\q_f^-$:
\begin{equation} \label{eq.path}
	\mathbf{\Phi}(\q^-, \q^+; s) = \mathbf{\Phi}^E(\q_E^+, \q_E^-; s) + \mathbf{\Phi}^f(\q_f^+, \q_f^-; s).
\end{equation}
This has a twofold effect: it allows to pointwise cancel the numerical errors introduced in the computation of the equilibrium and increases the accuracy of the solution out of the equilibrium.\\
In the context of MUSCL-Hancock strategy, a proper reconstruction $\q$ of the current numerical states $\Q$ is needed. In particular, following~\cite{castro2008well}, instead of reconstructing the current state, the reconstruction is performed on the fluctuations. Consequently, the reconstructed state is computed by summing the reconstructed fluctuation to the equilibrium, i.e.
\begin{equation} \label{eq.local_reconstruction_WB}
	\q_k(\x, t) = \Q_k^E(\x) + \mathbf{w}_k^f(\x, t) + \partial_t \Q_k^f (t - t^n) = \Q_k^E(\x) + \Q_k^f + \nabla \Q_k^f (\x - \x_k) + \partial_t \Q_k^f (t - t^n).
\end{equation} 
In~\eqref{eq.local_reconstruction_WB}, the spatial polynomial $\mathbf{w}_k^f$ is recovered through constraint~\eqref{eq.DeltaQk} but applied to the fluctuations $\Q_k^f = \Q_k^n - \Q_k^E$. Consequently, also the slope limiter~\eqref{eq.local_slopelimiter} is applied only to the discrete first order variation of the fluctuations, and the time coefficient $\partial_t \Q_k^f$ is given by~\eqref{eq.hyp} as
\begin{equation} \label{eq.Qt_f}
	\partial_t \Q_k = - \sum_{e_{kl} \in \mathcal{E}_k } | e_{kl} | \bigg( \F( \Q^E(\x_{kl}) + \mathbf{w}_k^f(\x_{kl}, t^n ) - \F( \Q^E(\x_{kl}) ) ) \cdot \n_{e_{kl}} \bigg) - \Big( \B(\Q_k^n) \cdot \nabla \Q_k - \B(\Q_k^E) \cdot \nabla \Q_k^E  \Big),  
\end{equation}
where from the flux and the nonconservative components along the reconstructed state we have subtracted the flux and the nonconservative components along the equilibrium. It follows that $\nabla \Q_k$ is the one used for recovering~\eqref{eq.local_reconstruction_WB}. Now, the scheme~\eqref{eq.FVWB} is effectively well balanced for any given equilibrium. As a matter of fact, at the equilibrium the initial condition $\Q_k^0 = \Q_k^E$ and the fluctuations are zero; 
consequently, the polynomial reconstruction of the state coincides with the equilibrium, i.e. $\q_k(\x) \equiv \Q_k^E(\x)$. Moreover, because of the smoothness of the equilibrium, the states from the two sides of the edge coincide, i.e. $\q_{e_{kl}}^- = \q_{e_{kl}}^+ = \Q_{e_{kl}}^E$. Thus the Rusanov flux~\eqref{eq.rusanov} reduces to be $\boldsymbol{\mathcal{F}}_{e_{kl}}(\q_{e_{kl}}^-, \q_{e_{kl}}^+) = \F(\Q_{e_{kl}}^E)$. Therefore, the scheme~\eqref{eq.FVWB} trivially becomes $\Q_k^{n+1} = \Q_k^n = \Q_k^E$. 

In some application, the equilibrium could be known only point-wise on the domain (e.g. example, in the case of recovering the fluid depth in geophysical applications). For this reason,  $\Q_k^E$ can be either the exact solution known for any $\x \in \Omega_k$ or its reconstruction from point-wise values.  
The well balanced approach presented in this subsection needs the local polynomial reconstruction~\eqref{eq.local_reconstruction_WB} of the state as the sum of the local reconstruction of the fluctuations around the given solution $\Q_k^E$ and the equilibrium itself. 
Successively the state is evolved via~\eqref{eq.FVWB} by subtracting the original PDE~\eqref{eq.hyp} from the scheme~\eqref{eq.FV} evaluated along~\eqref{eq.local_reconstruction_WB}. 
This process, thus, needs always to know the equilibrium as an input. 
In the next subsection, we propose some rearrangements in the original scheme~\eqref{eq.hyp} such that the knowledge of the exact equilibium profile is no more needed and the scheme is automatically well balanced for water at rest situations. 

\subsection{Well balanced scheme for water at rest equilibrium} \label{subsec.WBLR}

The well balanced scheme proposed in the previous subsection is valid for any smooth equilibrium and does not necessarily apply only to the problem under consideration of SW equations in GCC. Therefore, in this subsection, we attempt to generate a well balanced scheme in the same spirit of MUSCL-Hancock approach but specializing the solution $\Q^E$ to be water at rest~\eqref{eq.WR} (i.e. $\Q^E = \Q^{WR}$) and by exploiting the properties of the specific hyperbolic system~\eqref{eq.hyp} determined by fluxes~\eqref{eq.F} and matrices~\eqref{eq.B} collecting the nonconservative parts.\\

From~\eqref{eq.WR}, the water at rest solution reads $\Q^{WR} = [ h, 0, 0, b, \gamma_{11}, \gamma_{12}, \gamma_{22} ]^T$, with $h + b = \eta$ the constant free surface. {Fluxes~\eqref{eq.F} and block $\B_i^{met}$ of matrices of nonconservative components~\eqref{eq.B} strongly depend on the velocity. For this reason, when they are evaluated on $\Q^{WR}$, they vanish leading to a hyperbolic problem~\eqref{eq.hyp} that is flux-free and whose nonconservative components are only given by block $\B_i^\eta$.}\\

{Since the free surface $\eta$ is constant in space and time (condition~\eqref{eq.WRfs}), a second-order reconstruction exactly approximates it. Thus, it is convenient to consider the fluid depth $h$ defined by the difference of the free surface $\eta$ and the bathymetry $b$. In particular, let $b_k(\x)$ and $\eta_k(\x, t)$ be the local reconstructions of the bathymetry and the free surface, respectively. Reconstruction $b_k(\x)$ only depends on space because in our system it never evolves in time (so it can be reconstructed once and for all at the initial time). Moreover, concerning the bathymetry, we here consider the possibility to reconstruct it even though it is known a priori. For this reason, if the exact definition of $b$ is employed, $b_k^n(\x)$ will coincide with $b(\x)$ and $\nabla b_k$ will indicate the exact gradient $\nabla b(\x_k)$ evaluated in the cell-center $\x_k$.} Consequently, $h_k(\x, t) = \eta_k(\x, t) - b_k(\x)$ is the reconstruction of the fluid depth. At this point, smoothness hypothesis on the bathymetry can be relaxed, so we can also consider a bottom topography with jumps along certain points of the domain. From the definition of the local polynomial~\eqref{eq.local_reconstruction}, it follows that the modal coefficients for $h_k(\x, t)$ read
\begin{equation} \label{eq.hk}
	h_k^n = \eta_k^n - b_k^n, \quad \nabla h_k = \nabla \eta_k - \nabla b_k, \quad \partial_t h_k = \partial_t \eta_k.
\end{equation}

\medskip

In scheme~\eqref{eq.FV}, the nonconservative components $\B(\Q) \cdot \nabla \Q$ are approximated by $\B(\q_k^{n+\frac{1}{2}}) \cdot \nabla \Q_k$. Thus, through a polynomial reconstruction~\eqref{eq.hk}, we have
\begin{equation} \label{eq.ncp}
	\B_i^\eta(\q_k^{n + \frac{1}{2}}) \cdot \nabla \Q_k = \begin{bmatrix}
		0 \\
		\pm g h_k^{n + \frac{1}{2}} \gamma_{a b} (\partial_{x^i} h_k + \partial_{x^i} b_k) \\
		\mp g h_k^{n + \frac{1}{2}} \gamma_{c d} (\partial_{x^i} h_k + \partial_{x^i} b_k) \\
		\mathbf{0}_4
	\end{bmatrix},
\end{equation}
with $h_k^{n + \frac{1}{2}} = h_k(\x_k, t^{n + \frac{1}{2}})$, $a, b, c, d = 1, 2$ proper indexes for the metric, and $\partial_{x^i} h_k$ ($\partial_{x^i} b_k$) the $i$-th component of $\nabla h_k$ ($\nabla b_k$). Both nonzero entries of the right hand side of~\eqref{eq.ncp} are defined by $\partial_{x^i} h_k + \partial_{x^i} b_k$. By construction, at discrete level for these terms it holds
\begin{equation} \label{eq.detak}
	\partial_{x^i} h_k + \partial_{x^i} b_k = \partial_{x^i} \eta_k - \partial_{x^i} b_k + \partial_{x^i} b_k = \partial_{x^i} \eta_k.
\end{equation}
{It is now possible to remark the advantages in reconstructing the free surface $\eta$ instead of the fluid depth $h$. First, if the state is the water at rest equilibrium (i.e. $\eta$ is constant), then term~\eqref{eq.ncp} automatically vanishes because a constant function is discretized at machine precision $\epsilon_{mac}$ by a second-order polynomial reconstruction~\eqref{eq.local_reconstruction} (namely, $\nabla \eta_k = \mathcal{O}(\epsilon_{mac})$). Term~\eqref{eq.ncp} is thus automatically well balanced. Moreover, even though the bathymetry $b$ was discontinuous and therefore its local reconstruction $b_k(\x)$ had low resolution due to jumps, the reconstruction $h_k$ of the fluid depth $h$ relative to the free surface $\eta$ allows the contributions $\pm \partial_{x^i} b_k$ from bathymetry to be eliminated even at the \textit{discrete level}. Therefore, all numerical errors, eventually introduced by nonsmooth bottom topography through a (continuous) second-order polynomial~\eqref{eq.local_reconstruction}, are automatically cancelled in space and time. Thus, issues linked to discontinuities are only related to possible jumps of the free surface but they can never be caused by the bathymetry.}\\
\indent The same reasoning applies to the jumps at the interface~\eqref{eq.jump}. As a matter of fact, even after the integration of the  nonconservative matrices along the conservative path $\mathbf{\Psi}(\q^-, \q^+; \tau)$ of the states at left and at the right of an interface, the difference $\q^+ - \q^-$ with local reconstruction by~\eqref{eq.hk} leads to an equation of the form
\begin{equation} \label{eq.jumpWB}
	C \, (h^+ + b^+ - h^- - b^-) = C (\eta^+ - \eta^-),
\end{equation}
where $C$ is a term deriving from the integration of the matrix in~\eqref{eq.jump}. When the free surface is constant, the right hand side of~\eqref{eq.jumpWB} automatically vanishes.

\medskip

Finally, the last term to be analyzed in our scheme remains the discretization of the flux according to Rusanov flux~\eqref{eq.rusanov}. We have already noticed that in the case of zero kinetics, fluxes~\eqref{eq.F} automatically disappear. In formula~\eqref{eq.rusanov} then the components defined by $\F$ automatically vanish. However, for terms related to the difference $\q^- - \q^+$ of states to the left and right of the interface, they may not automatically cancel if a discontinuity in the bathymetry is present. \RIIcolor{In order to well balance this formula, we employ a slight different local stabilizer coefficient, as originally proposed in~\cite{castro2010some}. According to this approach, the used Rusanov solver reads
\begin{equation} \label{eq.rusanovWB}
	\boldsymbol{\mathcal{F}}_e (\q^-, \q^+) \cdot \n_e = \frac{1}{2} \left( \F(\q^+) + \F(\q^-) \right) \cdot \n_e - \frac{1}{2} s_{\max} \tilde{\mathbf{I}}_{mod} (\q^+ - \q^-),
\end{equation} }
{where $\tilde{\mathbf{I}}_{mod} \in \R^{m \times m}$ is the modified identity matrix introduced in~\eqref{eq.rusanov} with first row being $[1, 0, 0, 1, 0, 0, 0]$}. When the discrete flux~\eqref{eq.rusanovWB} is evaluated along the water at rest solution reconstruction, also the second part of the right hand side vanishes because the modified identity $\tilde{\mathbf{I}}$ makes the difference of states fall into the same situation as in formula~\eqref{eq.jumpWB}.

It is shown that through polynomial reconstruction~\eqref{eq.hk} and the use of a Rusanov flux~\eqref{eq.rusanovWB}, scheme~\eqref{eq.FV} is automatically well balanced for water at rest solution in the MUSCL-Hancock frame. The advantage of this approach consists on the fact that, when the water at rest solution is given as initial condition, the numerical system automatically preserves it. This approach is less invasive than the one presented in Sec.~\ref{subsec.WB}. Indeed, one just has to consider the reconstruction of the fluid depth depending on the reconstruction of the free surface. This option does not increase the number of operations to be performed compared to the not well balanced scheme. In fact, even for the latter, a local reconstruction of the conserved variables is required. In contrast, a classical well balanced approach would involve a reconstruction of the fluctuations and then the well-balancing of the scheme by subtracting the state evolution from the scheme at the equilibrium. 
Moreover, even outside equilibrium, such an approach allows states with discontinuous bathymetry to be considered without possible numerical errors arising from jumps affecting the recovered numerical solution. 
Finally, outside the equilibrium, the proposed minimally invasive well balanced numerical scheme is still second-order convergent. 

\section{Numerical results} \label{sec.test}
In this section, we exhibit the numerical examples (in both one and two dimensions) in support of previous sections. We start by showing how the scheme preserves the water at rest solution with smooth bathymetries and for different metrics, then we show numerical examples of both water at rest solutions with discontinuous bottom topographies and unsteady discontinuous states for the shallow water system~\eqref{eq.hyp}. In particular, we compare the effectiveness of the proposed well balanced scheme against the results obtained by employing a classical not well balanced scheme. Finally, some comparisons with results from the literature are presented.

\medskip

In all the simulations, the gravity constant $g$ is set equal to 9.81. Moreover, the $\CFL$ coefficient is always $0.9/d$, where $d=1,2$ is the space dimension of the specific test case.

The manifolds considered throughout the section refer to three cases of equipotential surfaces with respect to the gravitational field: the horizontal plane $\mathcal{M}_p$, the sphere $\mathcal{M}_R$ of radius $R > 0$ and the ellipse $\mathcal{M}_{K, \beta}$ of linear eccentricity $K > 0$ and constant surface level $\beta$: 
\begin{subequations} \label{eq.manifolds}
	\begin{alignat}{3}
		\mathcal{M}_p &= \R^2, \label{eq.manifolds_plane} \\
		\mathcal{M}_R &= \left\{ (x_1, x_2, x_3) \in \R^3: x_1^2 + x_2^2 + x_3^2 = R^2 \right\}, \label{eq.manifolds_sphere} \\
		\mathcal{M}_{K, \beta} &= \left\{ (x_1, x_2, x_3) \in \R^3: \frac{x_1^2}{\cosh^2 \beta} + \frac{x_2^2}{\cosh^2 \beta} + \frac{x_3^2}{\sinh^2 \beta} = K^2 \right\}, \label{eq.manifolds_ellipse}
	\end{alignat}
\end{subequations}
where $\x' = (x_1, x_2, x_3)$ are the Cartesian coordinates. Manifolds~\eqref{eq.manifolds} refer, therefore, to the Cartesian $\mathcal{H}_p$, the spherical $\mathcal{H}_R$ and the polar elliptic $\mathcal{H}_{K, \beta}$ maps
\begin{subequations} \label{eq.maps}
	\begin{alignat}{3}
		\mathcal{H}_p &= \left\{ \begin{array}{l}
			x_1 = x, \\
			x_2 = y, \\
			x_3 = z,
		\end{array} \right. & \x &= (x,y) \in [x_{\min}, x_{\max}] \times [y_{\min}, y_{\max}];	 \label{eq.maps_plane} \\
		\mathcal{H}_R &= \left\{ \begin{array}{l}
			x_1 = R \cos(\theta) \cos(\varphi), \\
			x_2 = R \sin(\theta) \cos(\varphi), \\
			x_3 = R \sin(\varphi),
		\end{array} \right. & \x &= (\theta,\varphi) \in [0, 2 \pi] \times \left[-\frac{\pi}{2}, \frac{\pi}{2} \right]; \label{eq.maps_sphere} \\
		\mathcal{H}_{K, \beta} &= \left\{ \begin{array}{l}
			x_1 = K \cosh(\beta) \cos(\theta) \cos(\varphi), \\
			x_2 = K \cosh(\beta) \sin(\theta) \cos(\varphi), \\
			x_3 = K \sinh(\beta) \sin(\varphi), 
		\end{array} \right. & \x &= (\theta,\varphi) \in [0, 2 \pi] \times \left[-\frac{\pi}{2}, \frac{\pi}{2} \right]; \label{eq.maps_ellipse}
	\end{alignat}
\end{subequations}
respectively. In maps~\eqref{eq.maps}, the covariant frame $\x = (x^1, x^2)$ is defined by the Cartesian coordinates $(x,y)$ in the Cartesian case, and by the longitudinal-latitudinal angles $(\theta, \varphi)$ for both the spherical and the elliptical contexts. The height $x_3$ in the frame of general covariant coordinates is given by the definition of the manifold. In particular, it specifies the normal direction to the manifold. In Cartesian coordinates, it coincides with the vertical direction parallel to $z$-axis; in spherical and elliptical coordinates, it is defined by the sinus of the latitude angle $\varphi$. For any map $\mathcal{H}$ in~\eqref{eq.maps}, let $J_i^j = \partial x_i / \partial x^j$ be the associated Jacobian matrix. Thus, the covariant tensor $\gamma_{ij}$ is  defined by the sum $J_i^k J^k_j$. In Tab.~\ref{tab.metrics} we report the metric tensor coefficients for the Cartesian, spherical and elliptical cases. These values are given as input data for the hyperbolic system~\eqref{eq.hyp}. Then, the numerical scheme (well balanced or not) discretizes their derivatives in order to properly recover the Christoffel symbols to account for the curvature of the manifold. The determinant of the metric tensor in Tab.~\ref{tab.metrics} is $\gamma = \gamma_{11} \gamma_{22} - \gamma_{12}^2$. For spherical and elliptical metrics, it vanishes at the north and south poles of the manifold. Consequently, these points cannot be considered as computational geometrical points.

In what follows, also one-dimensional test cases are presented. Let $\xi$ be the unique general covariant space coordinate for the test cases. We identify $\xi = x$ in the Cartesian frame and $\xi = \varphi$ (latitudinal angle) for both spherical and elliptical frames.

All quantitative analysis are performed through the $L^2$-norm
of the mismatch $\varepsilon_N$ between the recovered numerical solution $\Upsilon_N$ (with $\Upsilon_N = h, u^1, u^2$ dependent on the discretization $\mathcal{T}_N$) and the exact water at rest equilibrium $\Upsilon$ (with $\Upsilon_{ex} = h_{ex}, u^1_{ex}, u^2_{ex}$): $\varepsilon_N = \| \Upsilon - \Upsilon_{ex} \|_{L^2(\Omega)}$. Let $r_N$ be a characteristic mesh size depending on the partition $\mathcal{T}_N$. For two different discretizations $\mathcal{T}_{N_1}$ and $\mathcal{T}_{N_2}$, the order of convergence is computed as
\begin{equation} \label{eq.ooc}
	\mathcal{O}(L^2) = \frac{ \log(\varepsilon_{N_1} / \varepsilon_{N_2} ) }{ \log(r_{N_1} / r_{N_2}) }.
\end{equation}
We identify $r_N$ with the length cell $\Delta \xi$ in one dimension, and with the averaged incircle diameter $d_N$ of the polygonal elements in the tessellation $\mathcal{T}_N$, in two dimensions. 

\RIcolor{In order to make the analysis more convincing and challenging, a random perturbation is added to the equilibrium. In particular, the data collected in all the presented tables refer to double precision simulations. In some cases, explicitly specified, an additional analysis is also performed by working in single and quadruple precision.}

Finally, in order to speed up the simulations, a parallel implementations is performed both in 1D (exploiting the parallel MPI standard) and in 2D (thanks to the OpenMP interface).

\begin{table}[]
	\caption{Entries of the symmetric metric tensor in the Cartesian, spherical and elliptical frames. They derive from the Jacobian matrix $J_i^j$ associated to maps~\eqref{eq.maps}. The diagonal terms $\gamma_{11}$ and $\gamma_{22}$ with the extra-diagonal term $\gamma_{12} = \gamma_{21}$ are the only metric input data needed by the proposed hyperbolic model~\eqref{eq.hyp} to evolve the physical quantities and autonomously compute the curvature information of the manifold (and the metric information that we consider in the set of conserved variables). In both the spherical and the elliptical metrics, angle $\phi$ measures the latitude of the manifold. In the elliptical frame, $K$ is the linear eccentricity and $\beta$ is the constant surface level.}
	\label{tab.metrics}
	\centering
	\begin{tabular}{cccc} \toprule
		Metric     & $\gamma_{11}$                        & $\gamma_{12}=\gamma_{21}$ & $\gamma_{22}$                                                           \\ \midrule
		Cartesian  & 1                                    & 0             & 1                                                                       \\
		Spherical  & $R^2 \cos^2 \varphi$                 & 0             & $R^2$                                                                   \\
		Elliptical & $K^2 \cosh^2(\beta) \cos^2(\varphi)$ & 0             & $K^2(\cosh^2(\beta) \sin^2(\varphi) + \sinh^2(\beta) \cos^2(\varphi) )$ \\ \bottomrule
	\end{tabular}
\end{table}

\subsection{Preservation of water at rest equilibrium} \label{subsec.pres}
\begin{figure}
	\centering
	\subfigure[One-dimensional case of water at rest solution with free surface $\eta = 3$ (blue line) and Gaussian bump~\eqref{eq.bathy1} at the bottom (black line). The reference system is provided by the general covariant frame $(\xi, x^3)$, with $\xi = x$ in the Cartesian metric and $\xi = \varphi$ (latitudinal angle) for the spherical and elliptical metrics.]{
		\includegraphics[trim= 1 1 1 1,clip,width=0.45\linewidth]{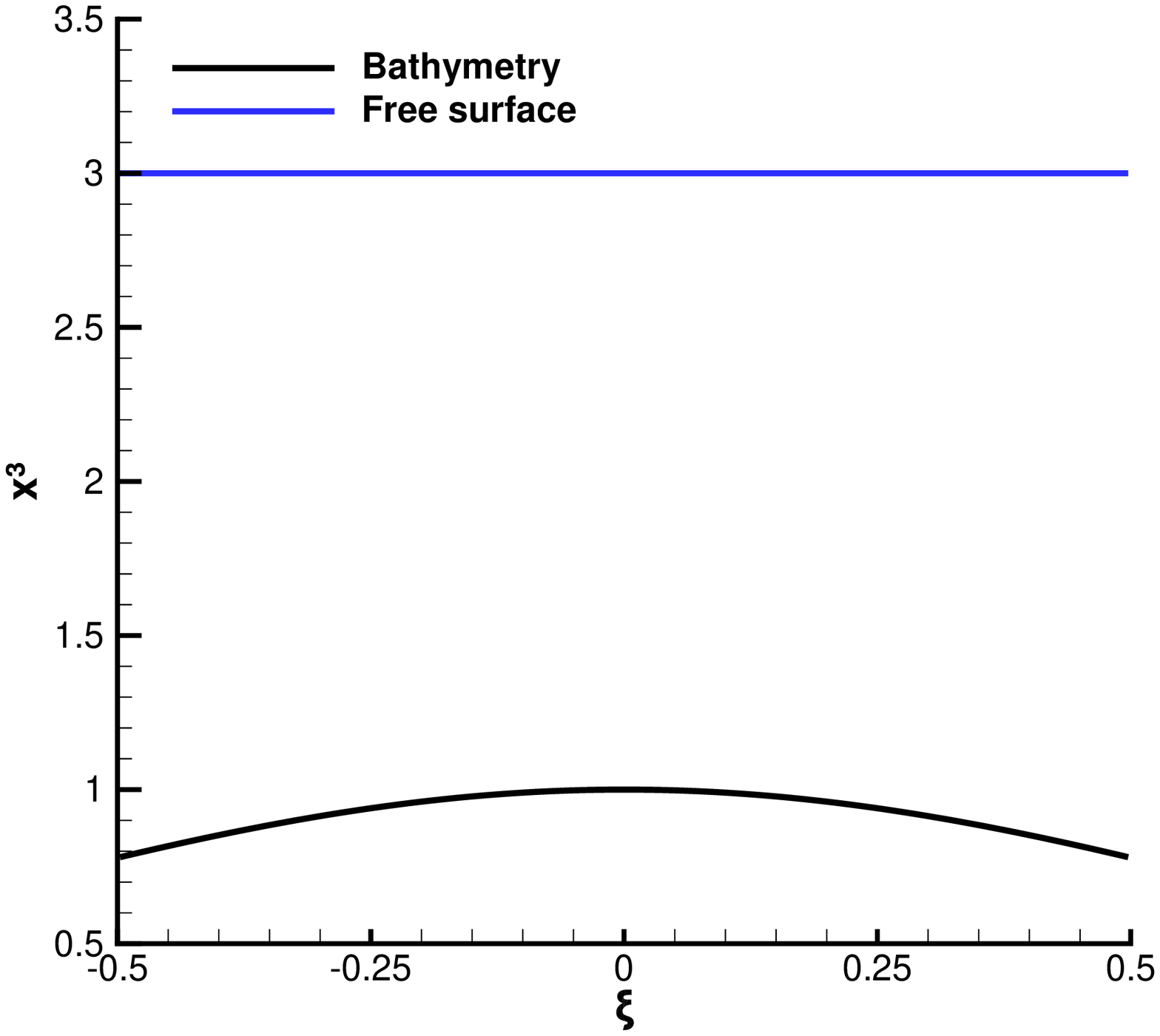}
		\label{fig.WR1D}
	}
	\hfill
	\subfigure[Two-dimensional case of water at rest solution with free surface $\eta = 3$ (surface on the top) and bump function ~\eqref{eq.bathy2} at the bottom (surface on the bottom). The reference system is provided by the general covariant frame $(x^1, x^2, x^3)$, with $(x^1, x^2) = (x, y)$ in Cartesian metric and $(x^1, x^2) = (\theta, \varphi)$ (longitudinal-latitudinal angles) for the spherical and elliptical metrics.]{
		\includegraphics[trim= 1 1 1 1,clip,width=0.45\linewidth]{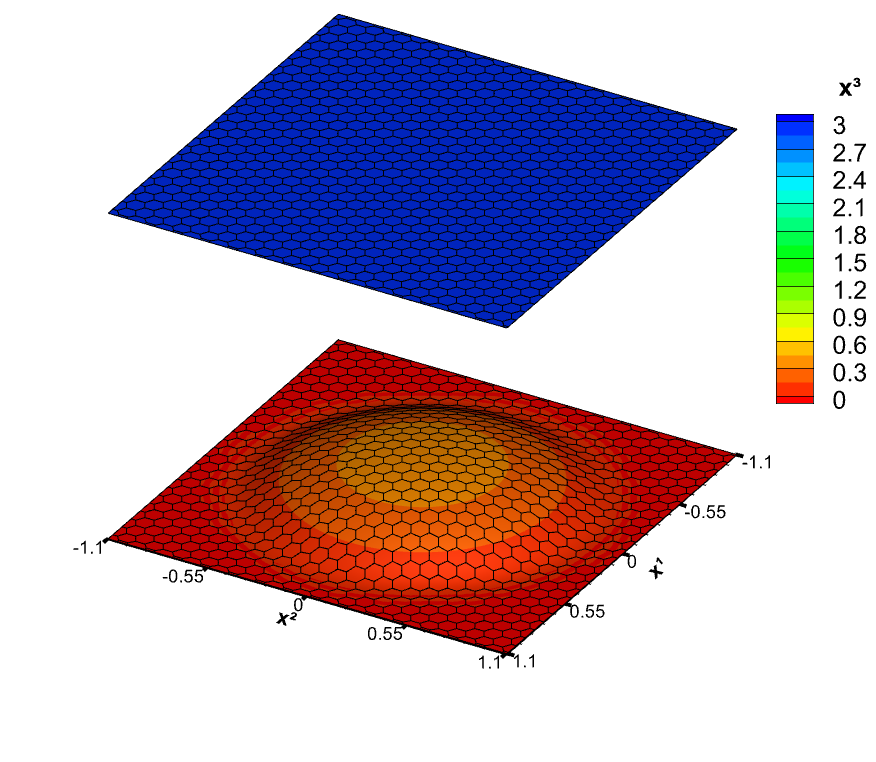}
		\label{fig.WR2D}
	}
	\caption{One-dimensional (a) and two-dimensional (b) water at rest equilibria in general covariant coordinates with free surface $\eta \equiv 3$ and bump functions~\eqref{eq.bathy} at the bottom.}
\end{figure}

The first set of numerical examples concern the ability of the scheme proposed in Sec.~\ref{subsec.WBLR} in preserving water at rest equilibria at machine precision. For the three metrics of Tab.~\ref{tab.metrics}, we consider a one-dimensional case and a two-dimensional case. For the one-dimensional case, the computational domain is $\Omega = [-0.5, 0.5]$ and, in the two-dimensional case, it turns into $\Omega = [-1.1, 1.1]^2$. The bathymetry is defined as
\begin{subequations} \label{eq.bathy}
	\begin{alignat}{2}
		b(\xi) &= e^{-\xi^2},	 \label{eq.bathy1} \\
		b(\x) = b(x^1, x^2) &= e^{-\frac{1}{1 - |\x|^2}} \chi_{ \{ |\x| \leq 1 \} } (\x) , \label{eq.bathy2} 
	\end{alignat} 
\end{subequations} 
for the one- and two-dimensional cases, respectively. 
In~\eqref{eq.bathy2}, $\chi_{ \{ |\x| \leq 1 \} } (\x)$ is the indicator function equal to 1 if $|\x| \leq 1$, 0 otherwise. In all cases, the free surface $\eta$ is set to be 3. For the spherical metric, the radius is $R = 1$. In elliptical coordinates, the linear eccentricity and the constant surface level measure $K = 1$ and $\beta = 2$, respectively. Fig.~\ref{fig.WR1D} and~\ref{fig.WR2D} show the one-dimensional and two-dimensional solutions, respectively. The pictures are plotted with respect to the reference system provided by the general covariant frame $(x^1, x^2, x^3)$. In particular, in both figures it is shown the bump on the bottom and the free surface of the water.

In Tab.~\ref{tab.WB_bGauss1D} we reported the $L^2$-errors on the physical variables of the 1D case. The values refer to cell length $\Delta \xi = 5.00$E-2. Also for large values of the final time $T$, it is evident that the water at rest equilibrium is preserved at machine-precision for all the considered metrics. 

For the two-dimensional case, Tab.~\ref{tab.WB_bGauss2D} resumes the $L^2$-errors for all physical variables, i.e. velocity components $(u^1, u^2)$ and fluid depth $h$. As long as the physical quantities are constant (as in the case of zero velocities), the scheme preserves the solution at machine precision even for long times. 
Concerning the fluid depth $h$, its variation in space is led by an exponential law defined by the bump function~\eqref{eq.bathy2} of the bathymetry. For this reason, it is no longer trivially constant throughout the computational domain $\Omega$. Considering the unstructured mesh of averaged size $d_N = 3.90E$-2 and despite machine precision is still preserved on large timescales, $L^2$-error on $h$ increases but never linearly, i.e. $\varepsilon_N^h = \| h - h_{ex} \|_{L^2(\Omega)} = o(t)$. 

This numerically proves the well-balancing property of the proposed scheme of Sec.~\ref{subsec.WBLR} for a smooth bathymetry. 

\begin{table}[]
	\caption{\RIcolor{$L^2$-error for 1D water at rest equilibrium through the proposed well balanced scheme of Sec.~\ref{subsec.WBLR} at varying of the used metrics (Cartesian, spherical and elliptical, see Tab.~\ref{tab.metrics}). The errors are measured for the fluid depth $h$ and the velocity $u$ for four time instances $t$. These values refer to a cell length $\Delta \xi = 5.00$E-2. Data are collected by working in single, double and quadruple precision. }  }
	\label{tab.WB_bGauss1D}
	\centering
	\begin{tabular}{cccccccc} \toprule
Precision & $t$ & \multicolumn{6}{c}{$L^2$-error} \\ \cmidrule(lr){3-8}
 &  & \multicolumn{2}{c}{Cartesian} & \multicolumn{2}{c}{Spherical} & \multicolumn{2}{c}{Elliptical} \\ \cmidrule(lr){3-4} \cmidrule(lr){5-6} \cmidrule(lr){7-8}
 &  & $h$ & $u$ & $h$ & $u$ & $h$ & $u$ \\ \midrule
\multirow{4}{*}{\rotatebox[origin=c]{90}{Single}} & 1 & 2.11E-7 & 1.84E-6 & 2.11E-7 & 1.84E-6 & 6.74E-8 & 2.22E-6 \\
 & 10 & 2.96E-7 & 2.36E-6 & 2.96E-7 & 2.36E-6 & 9.39E-8 & 4.34E-6 \\
 & 100 & 2.76E-7 & 2.69E-6 & 2.76E-7 & 2.69E-6 & 1.89E-7 & 6.01E-6 \\
 & 1000 & 2.96E-7 & 2.27E-6 & 2.96E-7 & 2.27E-6 & 1.88E-7 & 6.02E-6 \\ \midrule
\multirow{4}{*}{\rotatebox[origin=c]{90}{Double}} & 1 & 3.44E-13 & 1.63E-12 & 4.14E-13 & 1.93E-12 & 4.14E-13 & 1.83E-12 \\
 & 10 & 3.46E-13 & 1.63E-12 & 4.15E-13 & 1.93E-12 & 4.16E-13 & 1.77E-12 \\
 & 100 & 3.46E-13 & 1.63E-12 & 4.15E-13 & 1.93E-12 & 4.15E-13 & 1.78E-12 \\
 & 1000 & 3.46E-13 & 1.63E-12 & 4.15E-13 & 1.93E-12 & 4.15E-13 & 1.77E-12 \\ \midrule
\multirow{4}{*}{\rotatebox[origin=c]{90}{Quadruple}} & 1 & 4.99E-29 & 1.95E-29 & 4.99E-29 & 1.95E-29 & 4.97E-29 & 9.51E-30 \\
 & 10 & 5.11E-29 & 2.00E-29 & 5.11E-29 & 2.00E-29 & 5.06E-29 & 9.50E-30 \\
 & 100 & 4.91E-29 & 1.66E-29 & 4.91E-29 & 1.66E-29 & 5.01E-29 & 9.55E-30 \\
 & 1000 & 5.04E-29 & 1.64E-29 & 5.04E-29 & 1.64E-29 & 4.99E-29 & 1.30E-29 \\ \bottomrule
\end{tabular}
\end{table}

\begin{table}[]
	\caption{$L^2$-error for physical quantities $(h, u^1, u^2$ for 2D water at rest equilibrium through the proposed well balanced scheme of Sec.~\ref{subsec.WBLR} at varying of the used metric (Cartesian, spherical and elliptical, see Tab.~\ref{tab.metrics}). The considered bathymetry $b$ is defined in general covariant coordinates in~\eqref{eq.bathy2} for all metrics. Values in the table refer to four time instances $t$ and on a mesh whose characteristic size is $d_N = 3.90$E-2. }
	\label{tab.WB_bGauss2D}
	\centering
	\begin{tabular}{cccccccccc} \toprule
		$t$  & \multicolumn{9}{c}{$L^2$-error}                                                                  \\ \cmidrule(lr){2-10}
		& \multicolumn{3}{c}{Cartesian}  & \multicolumn{3}{c}{Spherical}  & \multicolumn{3}{c}{Elliptical} \\ \cmidrule(lr){2-4} \cmidrule(lr){5-7} \cmidrule(lr){8-10}
		& $h$      & $u^1$    & $u^2$    & $h$      & $u^1$    & $u^2$    & $h$      & $u^1$    & $u^2$    \\ \midrule
		1    & 3.83E-14 & 1.96E-14 & 1.53E-14 & 2.63E-14 & 5.50E-14 & 2.77E-14 & 2.10E-14 & 1.77E-14 & 9.22E-15 \\
		10   & 1.72E-13 & 2.32E-14 & 1.49E-14 & 3.90E-14 & 6.23E-14 & 3.46E-14 & 3.06E-14 & 1.21E-14 & 8.42E-15 \\
		100  & 1.49E-12 & 3.96E-14 & 3.02E-14 & 1.68E-13 & 1.18E-13 & 7.83E-14 & 6.21E-14 & 1.68E-14 & 1.18E-14 \\
		1000 & 1.50E-11 & 6.31E-14 & 5.22E-14 & 1.51E-12 & 1.32E-13 & 8.66E-14 & 4.23E-13 & 3.96E-14 & 2.55E-14	\\ \bottomrule
	\end{tabular}
\end{table}

\subsection{Discontinuous bottom topographies} \label{subsec:discbot}
In the previous subsection we considered continuous bathymetries for the performed analysis on the well balancing of the scheme. 
In this subsection, we relax the hypothesis on the continuity of the bottom topography. Throughout the entire subsection, only water at rest solutions will be considered. The goal is to numerically prove that the water at rest equilibrium is preserved due to the local polynomial reconstruction of the fluid depth $h$ as the difference between the reconstruction of the free surface $\eta$ and the bathymetry $b$.

\begin{figure}
	\centering
	\subfigure[Solution via non well balanced scheme.]{
		\includegraphics[trim= 2 2 2 2,clip,width=0.45\linewidth]{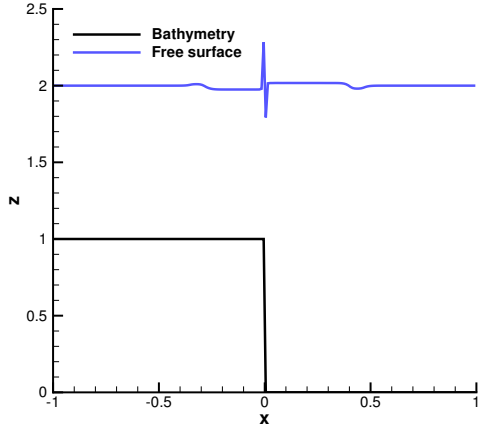}
	}
	\hfill
	\subfigure[Solution via well balanced scheme.]{
		\includegraphics[trim= 2 2 2 2,clip,width=0.45\linewidth]{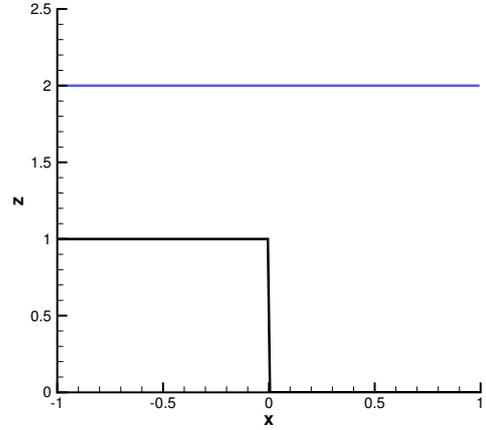}
	}
	\caption{Comparison of free surfaces for one-dimensional test case with Cartesian metric and jumping bathymetry along $x = 0$ obtained by both non well balanced scheme (a) and well balanced scheme (b). The final time and the cell length of mesh are $T = 0.1$ and $\Delta x = 1$E-2, respectively. }
	\label{fig:step1D_cart}
\end{figure}

Let us start with a one-dimensional test case in Cartesian metric whose domain is $\Omega = [-1,1]$ and bathymetry $b(x) = \chi_{\{ x \leq 0 \}}(x)$ defined by a step along $x = 0$. We set the free surface constantly equal to 2. Fig.~\ref{fig:step1D_cart} compares the free surface at final time $T = 0.1$  and cell length $\Delta x = 1$E-2 obtained by both the not well balanced scheme (on the left) and the well balanced scheme (on the right). The solution of the non well balanced scheme presents an nonphysical behavior around the jump in the bathymetry. In addition, due to this oscillation, two waves propagate backward and forward with respect to the jump.
On the contrary, the solution via well balanced scheme is preserved at machine precision. 

\begin{figure}
	\centering
	\subfigure[Cartesian]{
		\includegraphics[trim= 1 1 1 1,clip,width=0.45\linewidth]{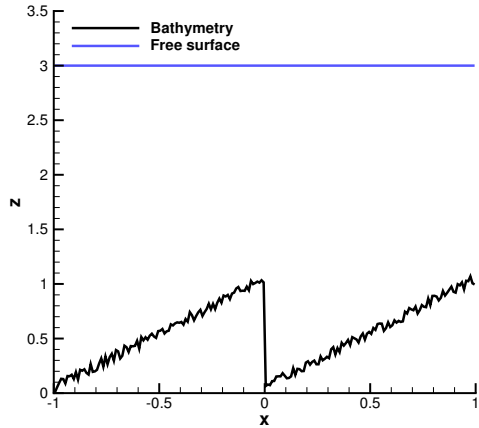}
	}
	\hfill
	\subfigure[Spherical]{
		\includegraphics[trim= 1 1 1 1,clip,width=0.45\linewidth]{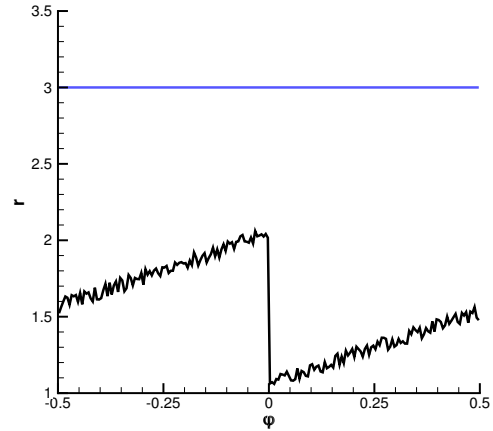}
	}
	\caption{Free surface for Cartesian (a) and spherical (b) test cases with noised jumping bathymetries~\eqref{eq.crazyb1D} via well balanced approach at final time $T = 1000$ and cell length $\Delta \xi = 1$E-2.}
	\label{fig.crazyb1D}
\end{figure}

To test the method, we considered two additional test cases in Cartesian and spherical metrics, respectively. For both cases, the original jumping bathymetries read
\begin{subequations} \label{eq.crazyb1D}
	\begin{alignat}{2}
		b(x) &= x + \chi_{\{ x \leq 0 \}}(x), & \quad x \in \Omega = [-1,1],	 \label{eq.crazyb1D_Cart} \\
		b(\varphi) &= \sin(\varphi) + 2\chi_{\{ \varphi \leq 0 \}}(\varphi) + \chi_{\{ \varphi > 0 \}}(\varphi), & \quad \varphi \in \Omega = [-0.5,0.5], \label{eq.crazyb1D_Pol} 
	\end{alignat} 
\end{subequations}
for the Cartesian metric~\eqref{eq.crazyb1D_Cart} and the spherical metric~\eqref{eq.crazyb1D_Pol}, respectively. Successively, to the expressions in~\eqref{eq.crazyb1D}, a white noise is added in order to simulate a discontinuous bottom topography in any cell (see black lines in Fig.~\ref{fig.crazyb1D}). For both cases, the free surface $\eta$ is equal to 3 and the cell length $\Delta \xi = 1$E-2.  
Because of the discontinuities, a non well balanced scheme is not able to catch the equilibrium and the scheme becomes unstable after few time instances. 
On the contrary, the proposed well balanced approach preserves the equilibria, as showed in Tab.~\ref{tab.crazyb1D} by the $L^2$-errors on the physical quantities at final time $T = 1000$.

\begin{figure}
	\centering
	\subfigure[Solution via non well balanced scheme. $T = 0.1$.]{
		\includegraphics[scale=.07]{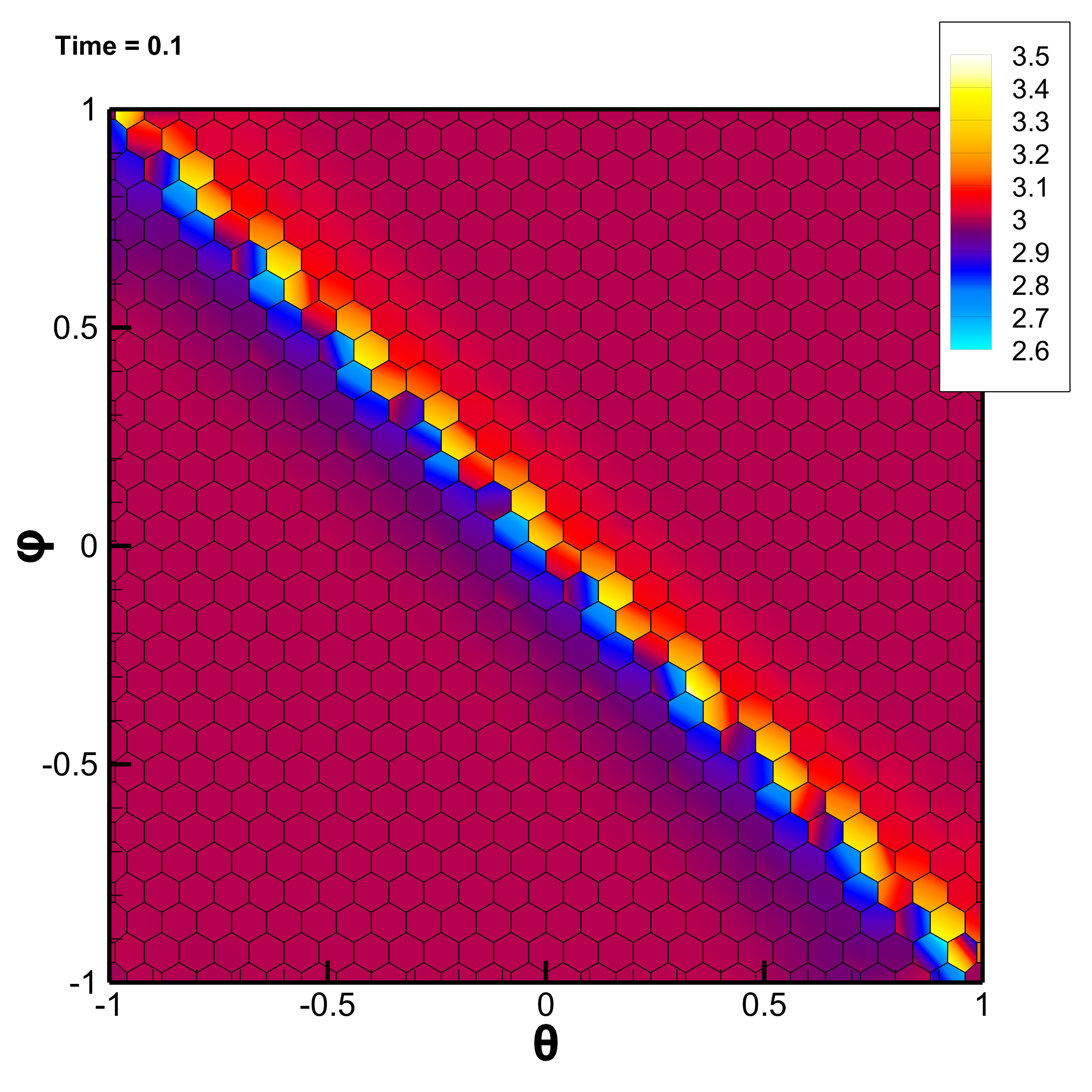}
	}
	\subfigure[Solution via non well balanced scheme on a finer mesh. $T = 0.1$.]{
		\includegraphics[scale=.07]{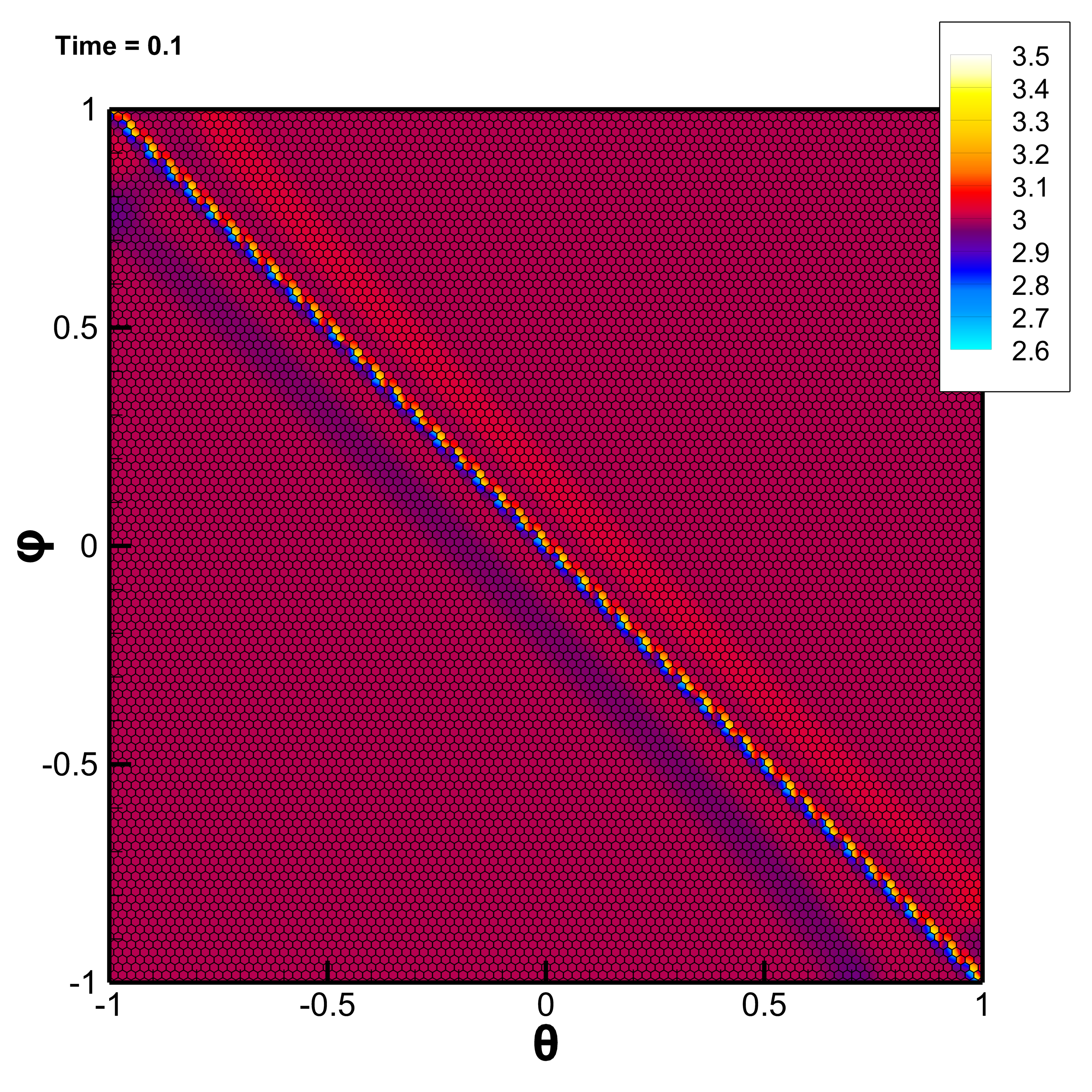}
	}
	\subfigure[Solution via well balanced scheme. $T = 10$.]{
		\includegraphics[scale=.07]{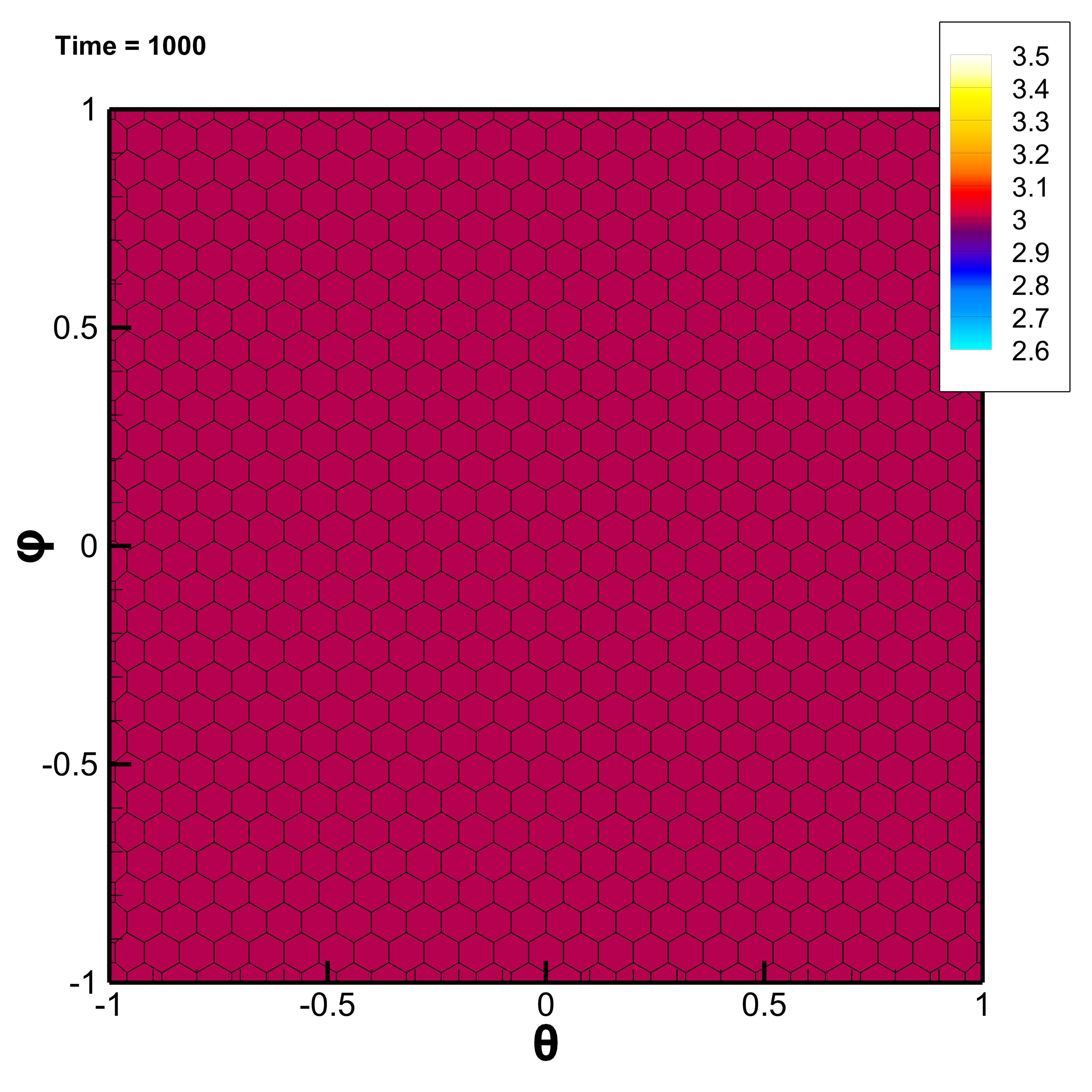} \label{fig.ellWB}
	}
	\caption{Comparison of the free surface via the non well balanced (a) and (b) well balanced (c) schemes for water at rest solution with free surface $\eta = 3$ and discontinuous bathymetry~\eqref{eq.ellbat} in elliptical coordinates. Figures (a) and (c) refer to an averaged mesh size $d_N = 2.85$E-2, figure (c) has finer mesh with $d_N = 7.39E$-3. The final time $T$ is 0.1 for figures (a) and (b) and 1000 for figure (c).}
	\label{fig.ell}
\end{figure}

We perform a similar analysis also in two dimensions in elliptical metric. The discontinuous bathymetry in covariant coordinates reads
\begin{equation} \label{eq.ellbat}
	b(\theta, \varphi) = \theta + \varphi + \chi_{\{\theta + \varphi \geq 0\}}(\theta, \varphi)
\end{equation}
that presents a discontinuity along the straight line $\varphi = - \theta$. The free surface is set to be equal to 3. Tab.~\ref{tab.ellconv} records the $L^2$-errors for the three physical variables and up to $T = 1000$ the solution is preserved at machine precision as for water at rest solution fully continuous. 
If the not well balanced scheme is employed, the resting phenomenon is not caught neither by refining the mesh (see Fig.~\ref{fig.ell}). Indeed, nonphysical oscillations appear around the discontinuity in the bathymetry that affect the precision also far from the discontinuity line.
On the contrary, even with a coarse mesh (Fig.~\ref{fig.ellWB}) the well balanced approach preserves the equilibrium on long timescales.

\begin{table}[]
	\caption{\RIcolor{$L^2$-error for one-dimensional test cases with noised bathymetries~\eqref{eq.crazyb1D} at time $T = 1000$. Data are collected by working in single, double and quadruple precision.}}
	\label{tab.crazyb1D}
	\centering
	\begin{tabular}{ccccc} \toprule
Precision & \multicolumn{4}{c}{$L^2$-error} \\ \cmidrule(lr){2-5}
 & \multicolumn{2}{c}{Cartesian} & \multicolumn{2}{c}{Spherical} \\ \cmidrule(lr){2-3} \cmidrule(lr){4-5}
 & $h$ & $u$ & $h$ & $u$ \\ \midrule
Single & 2.50E-7 & 3.03E-6 & 2.38E-8 & 1.74E-6 \\
Double & 7.06E-14 & 2.64E-13 & 9.21E-14 & 4.04E-13 \\
Quadruple & 5.80E-34 & 6.59E-33 & 4.87E-34 & 3.76E-33 \\ \bottomrule
\end{tabular}
\end{table}

\begin{table}[]
	\caption{$L^2$-error for physical quantities $(h, u^1, u^2)$ for 2D water at rest equilibrium through the proposed well balanced scheme with elliptical metric, free surface $\eta = 3$ and discontinuous bathymetry~\eqref{eq.ellbat}. The errors vary with respect to time $t$. The used mesh has an averaged mesh size $d_N = 2.85E$-2.}
	\label{tab.ellconv}
	\centering
	\begin{tabular}{cccc} \toprule
		$t$  & \multicolumn{3}{c}{$L^2$-error} \\ \cmidrule(lr){2-4}
		& $h$     & $u^1$     & $u^2$     \\ \midrule
		1    & 3.39E-14      & 4.74E-14       & 5.09E-14       \\
		10   & 3.22E-14      & 4.75E-14       & 4.97E-14       \\
		100  & 3.22E-14      & 5.11E-14       & 4.98E-14       \\
		1000 & 3.20E-14      & 4.40E-14       & 5.04E-14      \\ \bottomrule
	\end{tabular}
\end{table}

\subsection{Discontinuous solutions out of equilibrium} \label{subsec.discsol}
The trick of approximating the local reconstruction of the fluid depth $h$ by subtracting the bathymetry $b$ (eventually badly approximated due to the discontinuities) to the free surface $\eta$ not only allows to well balance the scheme but also to avoid nonphysical oscillations for states out of the equilibrium.

\begin{figure}
	\centering
	\subfigure[Solution via non well balanced scheme.]{
		\includegraphics[trim= 0 0 0 0,clip,width=0.45\linewidth]{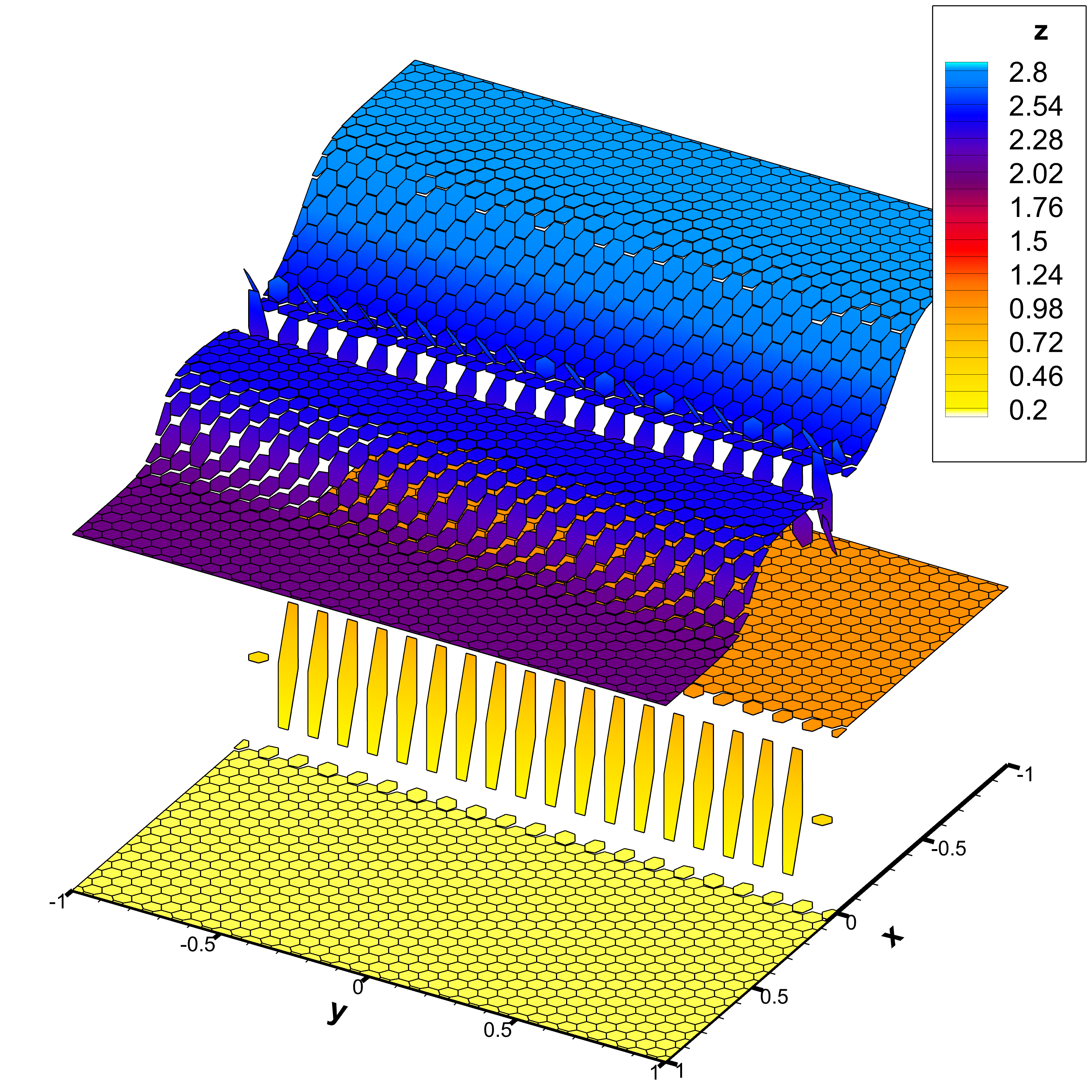} \label{fig.step2D_Cart_noWB}
	}
	\hfill
	\subfigure[Solution via well balanced scheme.]{
		\includegraphics[trim= 0 0 0 0,clip,width=0.45\linewidth]{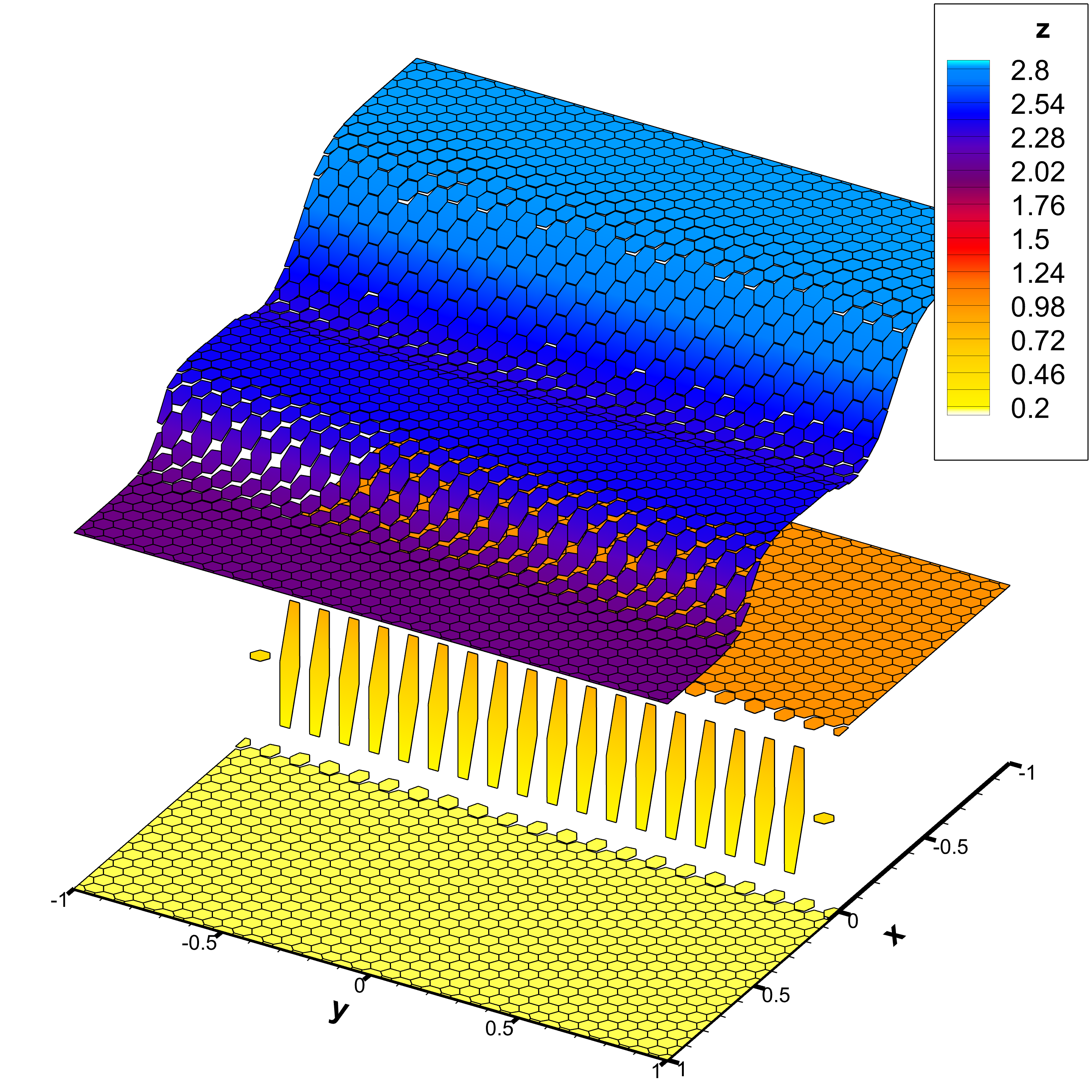} \label{fig.step2D_Cart_WB}
	}
	\caption{Comparison of free surface for two-dimensional test case with Cartesian metric and jumping bathymetry along line $x = 0$ (bottom surface) obtained by both non well balanced scheme (a) and well balanced scheme (b). The final time and the averaged characteristic length of the mesh are $T = 0.1$ and $d_N = 2.08$E-2, respectively.}
	\label{fig.step2D_Cart}
\end{figure}

The considered test case has Cartesian metric with a step at the bottom (i.e. $b(x,y) = \chi_{\{ x \leq 0 \}}(x,y)$). The physical variables respond to a Riemann problem with initial conditions defined by zero velocity and a discontinuous fluid depth such that the free surface is $\eta(x,y,t = 0) = 2 + \chi_{\{ x \leq 0 \}}(x,y)$. In Fig.~\ref{fig.step2D_Cart} we compare the simulations at final time $T = 0.1$ for the free surface (top surface) via the non well balanced (on the left) and well balanced (on the right) schemes on a mesh whose averaged characteristic length is $d_N = 2.08$E-2. 
The solution obtained through the non well balanced approach presents nonphysical oscillations along the discontinuity line $x = 0$ in bathymetry. The oscillations disappear when the proposed well balanced approach is adopted because in the scheme the contributions of the reconstructed bathymetry cancel at discrete level, as remarked in Sec.~\ref{subsec.WBLR}. 

\subsection{Order of convergence} \label{subsec.order}
\begin{figure}
	\centering
	\subfigure[Fluid depth $h$, bathymetry $b$ and relative free surface $\eta$]{
		\includegraphics[scale=.35]{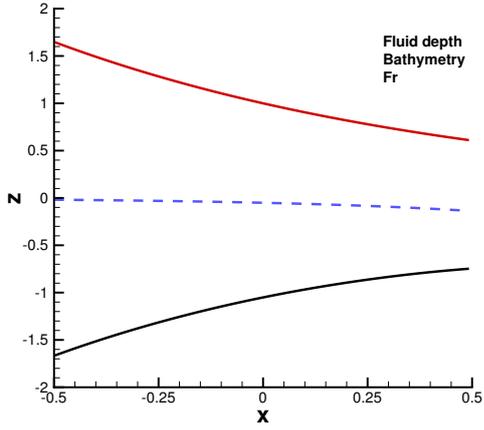}
	}
	\hfill
	\subfigure[Velocity $u$]{
		\includegraphics[scale=.35]{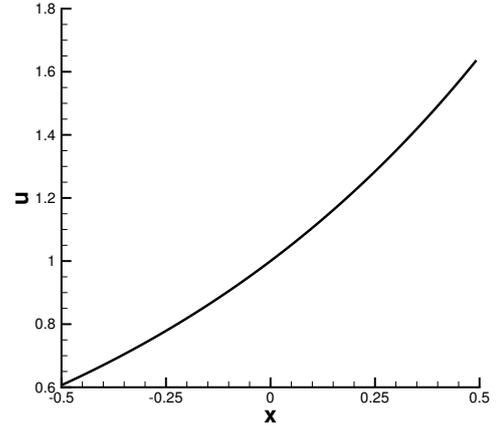}
	}
	\caption{Fluid depth $h$, bathymetry $b$ and relative free surface $\eta$ (a) and velocity $u$ (b) of the steady solution~\eqref{eq.particularsteadysolution} used for performing the order of convergence of the proposed well balanced scheme.}
	\label{fig.steady_solution}
\end{figure}

So far, we have verified numerically that the well balanced method we propose preserves water at rest equilibria and allows discontinuous solutions (even out of this equilibrium) not to exhibit nonphysical oscillations due to discontinuities in bottom topography. 
Now we take care of verifying that the well balanced method also converges with order two.
The order of convergence is computed on a particular one-dimensional class of steady solution for the original problem~\eqref{eq.baldauf} in Cartesian coordinates. This class is defined by any choice of fluid depth, velocity and bathymetry $(h, u, b)$ that fulfill 
\begin{equation} \label{eq.steadysolution}
	\begin{aligned}
		& \partial_t h = 0, \quad \partial_t u = 0, \\
		& h u \equiv 0, \quad h \neq 0, \quad \forall x \in \Omega, \\
		& \partial_x b = \left( \frac{u^2}{g h} -1 \right) \partial_x h.
	\end{aligned}
\end{equation}
The water at rest class of equilibria is a particular subset of solutions~\eqref{eq.steadysolution}. The chosen solution in~\eqref{eq.steadysolution} for  the order of convergence reads
\begin{equation} \label{eq.particularsteadysolution}
	h(x) = e^{-x}, \quad u(x) = e^x, \quad b(x) = -\frac{e^{2x}}{2 g} - e^{-x}.
\end{equation} 
Fig.~\ref{fig.steady_solution} depicts the bathymetry, the fluid depth and the relative free surface (on the left) and the velocity (on the right) of the steady solution~\eqref{eq.particularsteadysolution}. Tab.~\ref{tab.ooc} sums up the order of convergence for the fluid depth $h$ and the velocity $u$. The scheme is converging with rate 2 as expected from theory.

\begin{table}[]
	\caption{Table of convergence for equilibrium~\eqref{eq.particularsteadysolution}. The second and fourth columns report the $L^2$-error for fluid depth $h$ and velocity $u$ at varying of the cell length $\Delta x$, respectively. Columns 3 and 5 measure the respective convergence rates.}
	\label{tab.ooc}
	\centering
	\begin{tabular}{ccccc} \toprule
		& \multicolumn{2}{c}{$h$}        & \multicolumn{2}{c}{$u$}        \\ \cmidrule(lr){2-3} \cmidrule(lr){4-5}
		$\Delta x$ & $L^2$-err & $\mathcal{O}(L^2)$ & $L^2$-err & $\mathcal{O}(L^2)$ \\ \midrule
		2.00E-2    & 1.1698E-5 & -                  & 1.8417E-5 & -                  \\
		1.00E-2    & 2.9651E-6 & 1.98               & 4.6648E-6 & 1.98               \\
		5.00E-3    & 7.4627E-7 & 1.99               & 1.1730E-6 & 1.99               \\
		3.33E-3    & 3.3241E-7 & 2.00               & 5.2235E-7 & 2.00               \\
		2.50E-2    & 1.8718E-7 & 2.00               & 2.9411E-7 & 2.00              \\ \bottomrule
	\end{tabular}
\end{table}

\subsection{Metric-adaptive property of the model} \label{subsec.metric}
We close this section by numerically analyzing one main property of the hyperbolic model~\eqref{eq.hyp} derived from the model of~\eqref{eq.baldauf}. Writing it in this hyperbolic form, with the metric tensor taken in the set of conserved variables, at the same time allows it to be easily discretized with a finite volume scheme but, more importantly, to avoid having to explicitly compute the curvature of the manifold (hence Christoffel symbols). 
In particular, the model needs as its initial input data the components of the metric tensor, and then the system autonomously approximates the geometric characteristics of the manifold. 
In this subsection we want to compare numerical results of this model with results derived by other models built \textit{ad hoc} for a specific manifold.

The first comparison is performed in Cartesian metric. In this context the manifold trivially degenerates to the plane $\R^2$. Consequently, the shallow water system simply reads
\begin{equation} \label{eq.SWECartesian}
	\begin{aligned}
		\partial_t h + \partial_{x_j} m_j &= 0, \\
		\partial_t m_i + \partial_{x_j} \left( \frac{m_i m_j}{h} + \frac{1}{2} g h^2 \right) &= - g h \partial_{x_j} b,
	\end{aligned}
\end{equation}
where $(x_1, x_2) = (x,y)$ are the classical Cartesian coordinates. 
System~\eqref{eq.SWECartesian} is already written in hyperbolic form whose conserved variables are $\Q = [h, m_1, m_2]^T$. 
It is possible to prove that the associated Jacobian matrix has a spectrum that coincides with~\eqref{eq.eigenvalues} of our model. In particular, the characteristic velocity $c$ is equal to $\sqrt{g h}$. We noticed in Remarks~\ref{rem.cartesian2} and~\ref{rem.cartesianWR} that this theoretically implies that the numerical solutions of systems~\eqref{eq.SWECartesian} and~\eqref{eq.hyp} (with metric $\gamma_{ij} = \delta_{ij}$) have to coincide at machine precision. 

\begin{figure}
	\centering
	\subfigure[Fluid depth $h$]{
		\includegraphics[trim= 1 1 1 1,clip,width=0.45\linewidth]{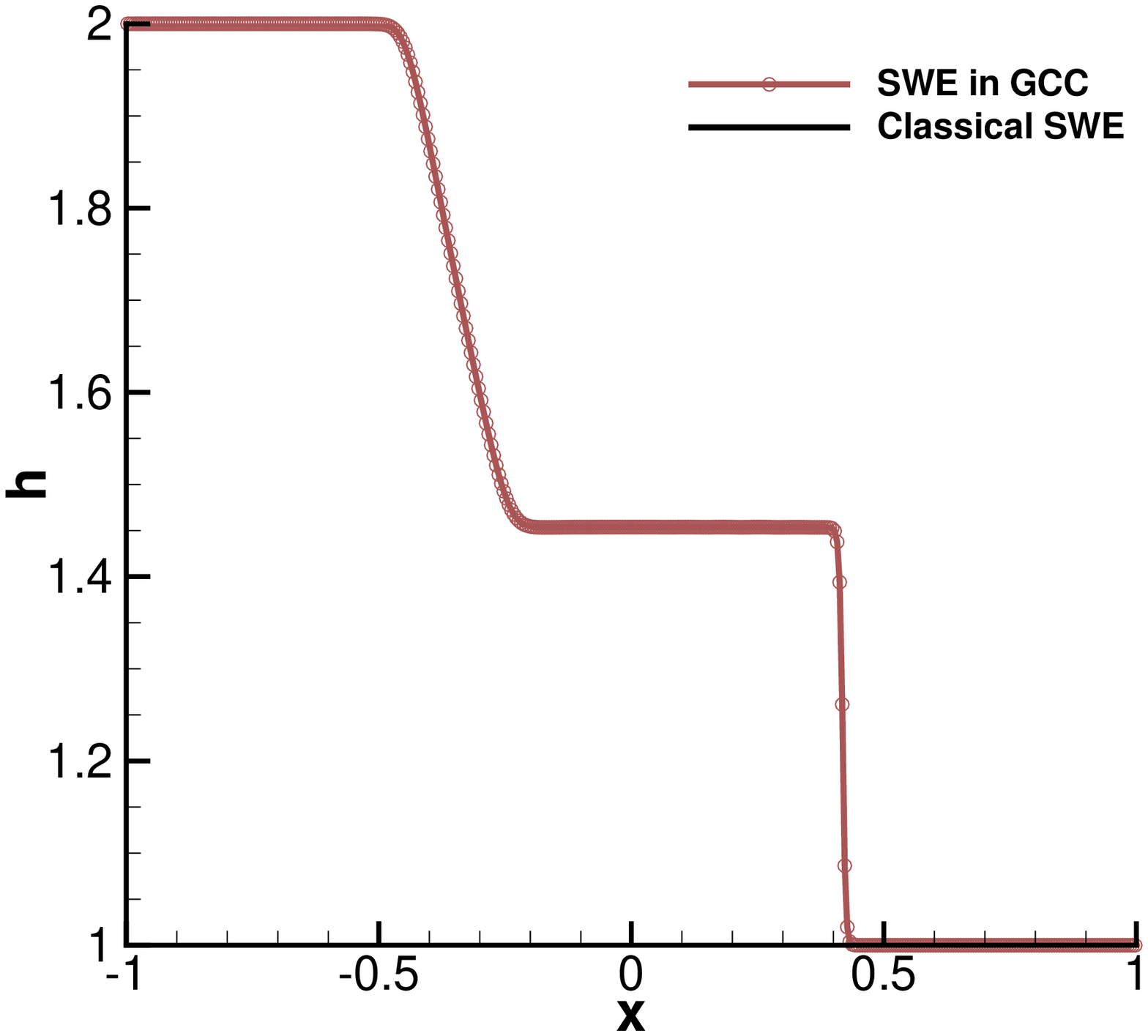}
	}
	\subfigure[Velocity $u$]{
		\includegraphics[trim= 1 1 1 1,clip,width=0.45\linewidth]{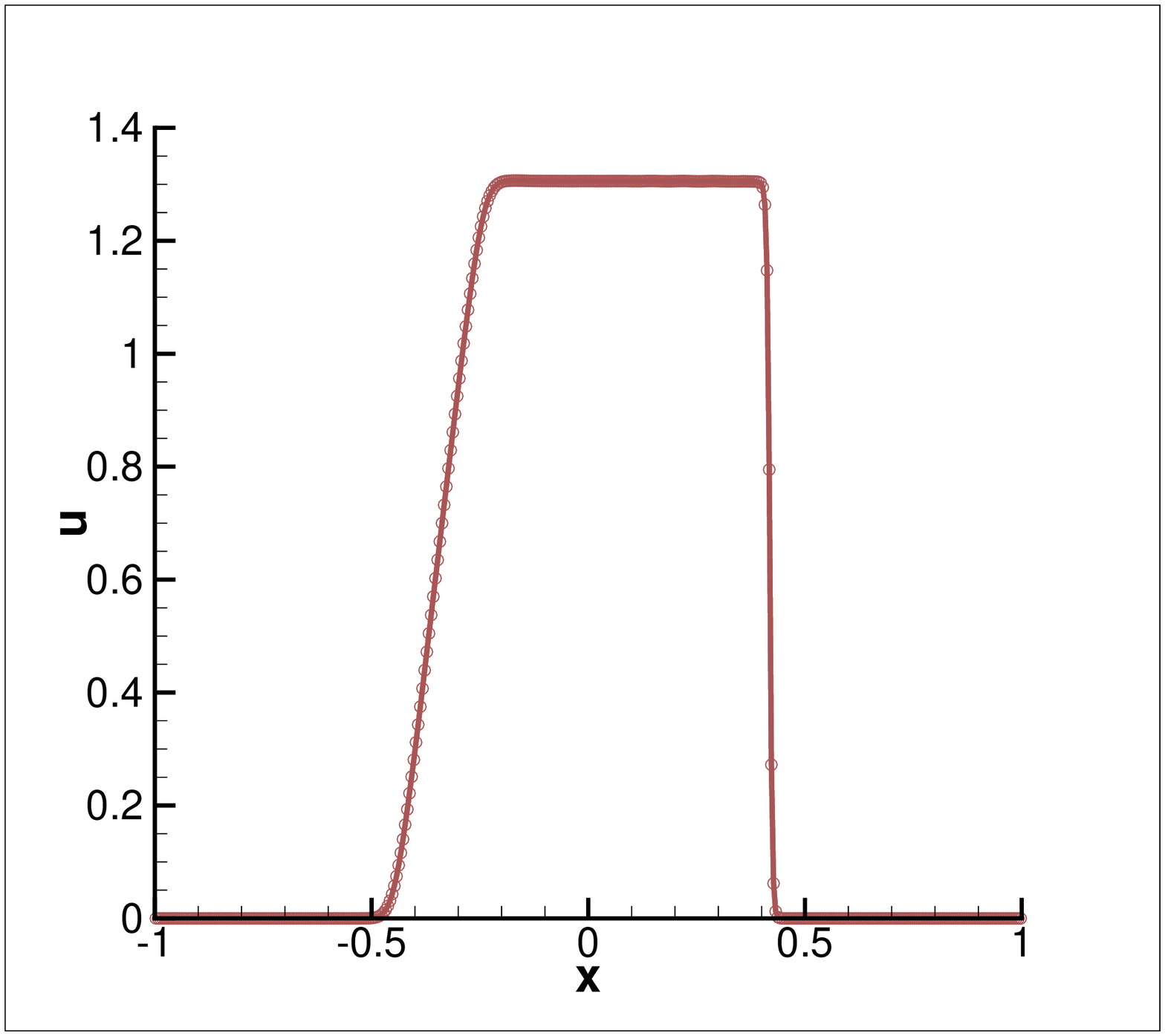}
	}
	\caption{Comparison of a Riemann problem for solutions from the proposed hyperbolic model~\eqref{eq.hyp} in Cartesian metric and from classical model~\eqref{eq.SWECartesian} for physical quantities $h$ (a) and $u$ (b). The final time and the cell length of the mesh are $T = 0.1$ and $\Delta x = 1$E-2, respectively.}
	\label{fig.SWECartesian}
\end{figure}

Let $\Omega = [-1,1]$ be the one-dimensional computational domain. The solved problem is a Riemann problem with a flat bathymetry $b \equiv 0$, zero initial velocity and an initial fluid depth $h(x, t = 0) = 2 + \chi_{\{ x \leq 0 \}}(x)$ presenting a discontinuity at the origin $x = 0$. In Fig.~\ref{fig.SWECartesian} there are the evolved curves of fluid depth $h$ and velocity $u$ at final time $T = 0.1$ for a mesh whose cell length $\Delta x = 1$E-2. We overlap the solutions from the general model~\eqref{eq.hyp} in Cartesian metric and the classical model~\eqref{eq.SWECartesian}. At a glance there is no difference. As a matter of fact the mismatch of the solutions of the two models in $L^2$-norm is of the order of $10^{-15}$.
This numerically proves that the proposed hyperbolic model~\eqref{eq.hyp} automatically degenerates to the classical SW model~\eqref{eq.SWECartesian} when the metric is set to be Cartesian.

\begin{figure}
	\centering
	\subfigure[Solution from model~\eqref{eq.hyp} with spherical metric.]{
		\includegraphics[trim= 1 1 1 5,clip,width=0.45\linewidth]{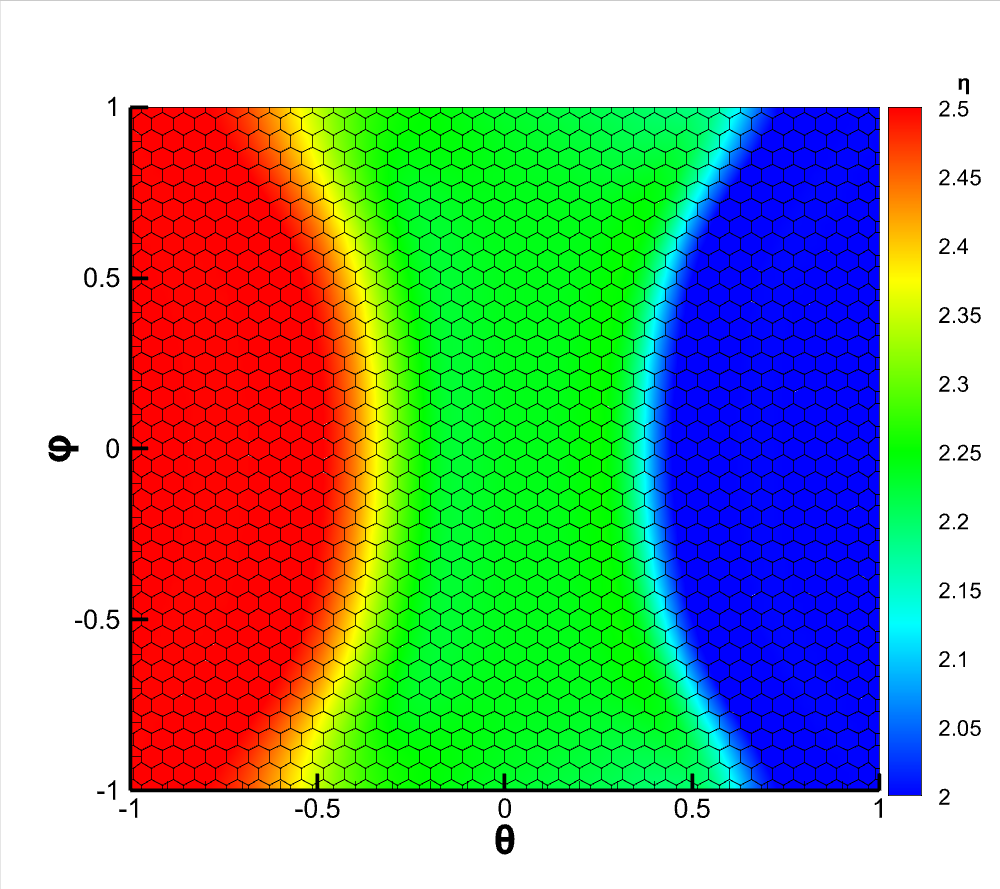}
	}
	\hfill
	\subfigure[Solution from model~\eqref{eq.castro}~\cite{castro2017well}.]{
		\includegraphics[trim= 1 1 1 5,clip,width=0.45\linewidth]{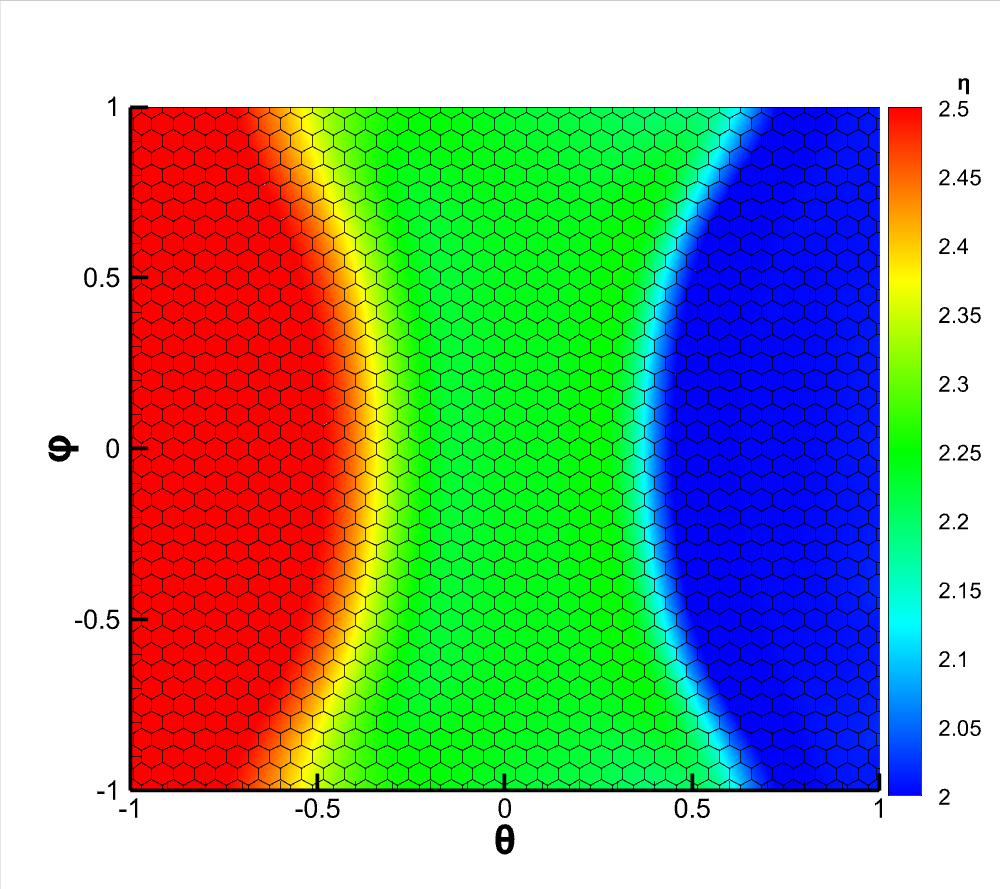}
	}
	\caption{Comparison of solutions of a Riemann problem in spherical coordinates from the proposed model~\eqref{eq.hyp} and model~\eqref{eq.castro} presented by Castro \textit{et. al} in~\cite{castro2017well}. For both cases, the final time and the averaged mesh size are $T = 0.1$ and $d_N = 5E$-2 . }
	\label{fig.SWEPolar}
\end{figure}

\begin{figure}
	\centering
	\subfigure[Solution from model~\eqref{eq.hyp} with spherical metric.]{
		\includegraphics[width=.3\linewidth]{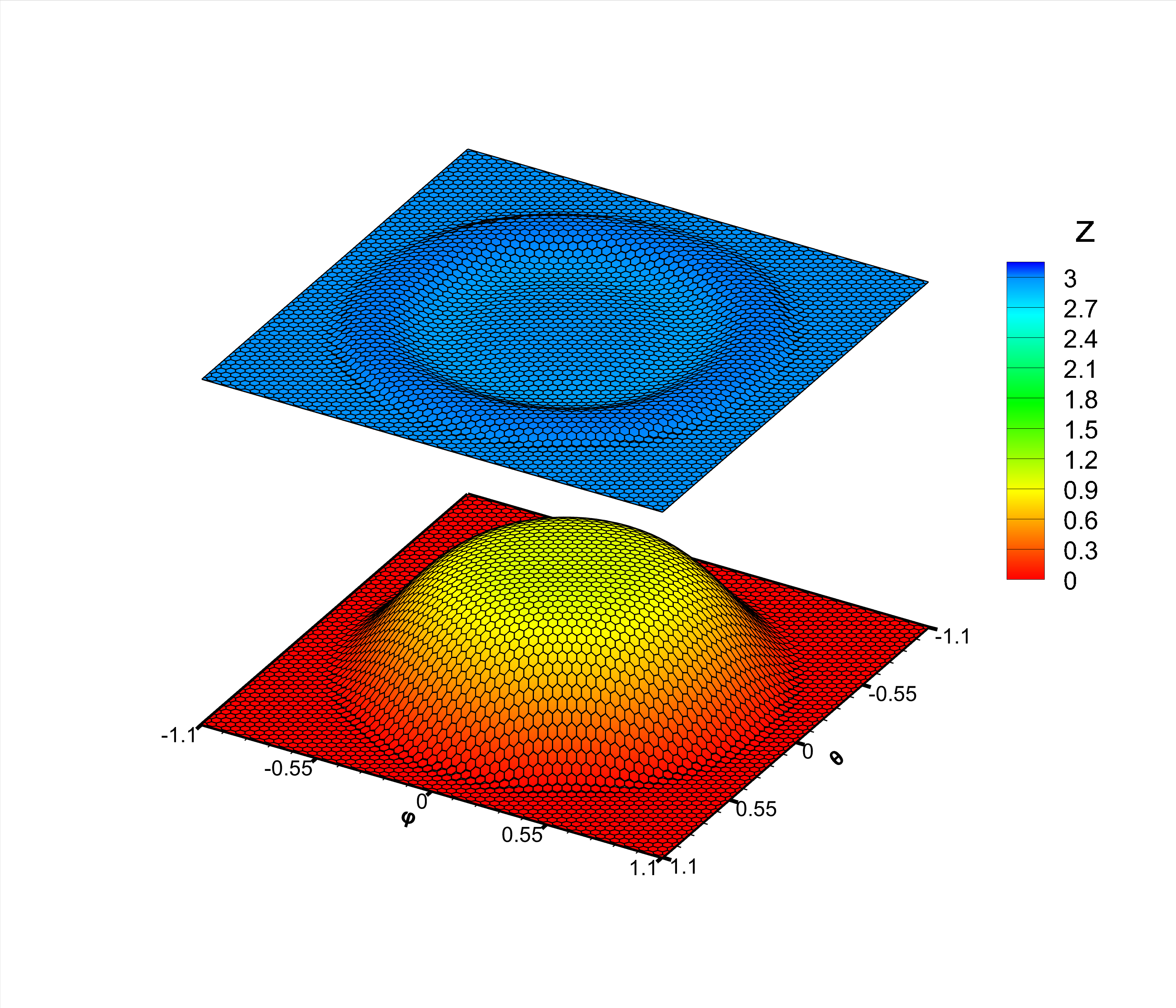}
	}
	\hfill
	\subfigure[Solution from model~\eqref{eq.castro}~\cite{castro2017well}.]{
		\includegraphics[width=.3\linewidth]{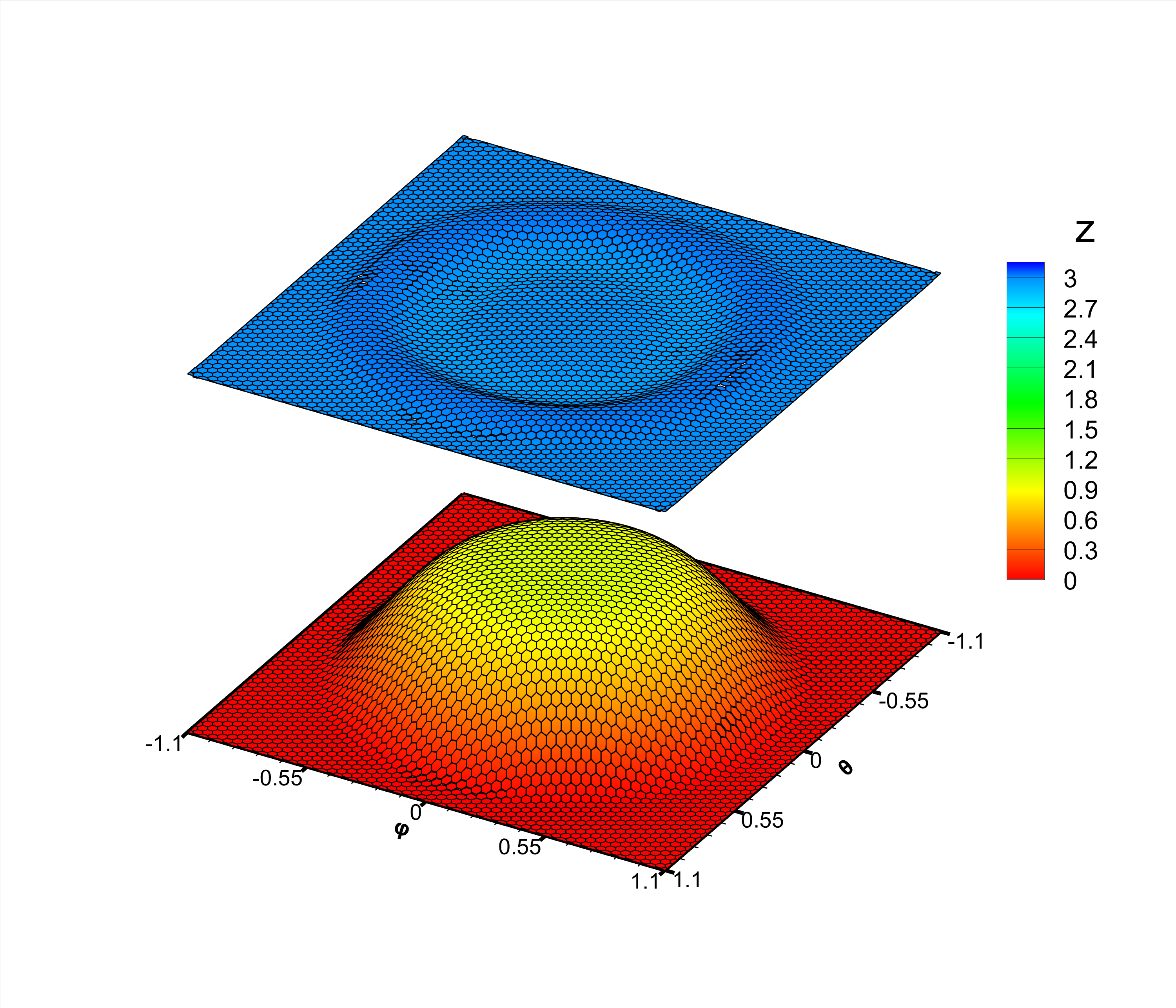}
	}
	\hfill
	\subfigure[Free surface profile comparison for $\theta \in (-1,1)$ and $\varphi = 0$.]{
		\includegraphics[trim = 8 6 6 6 , clip, width=.25\linewidth]{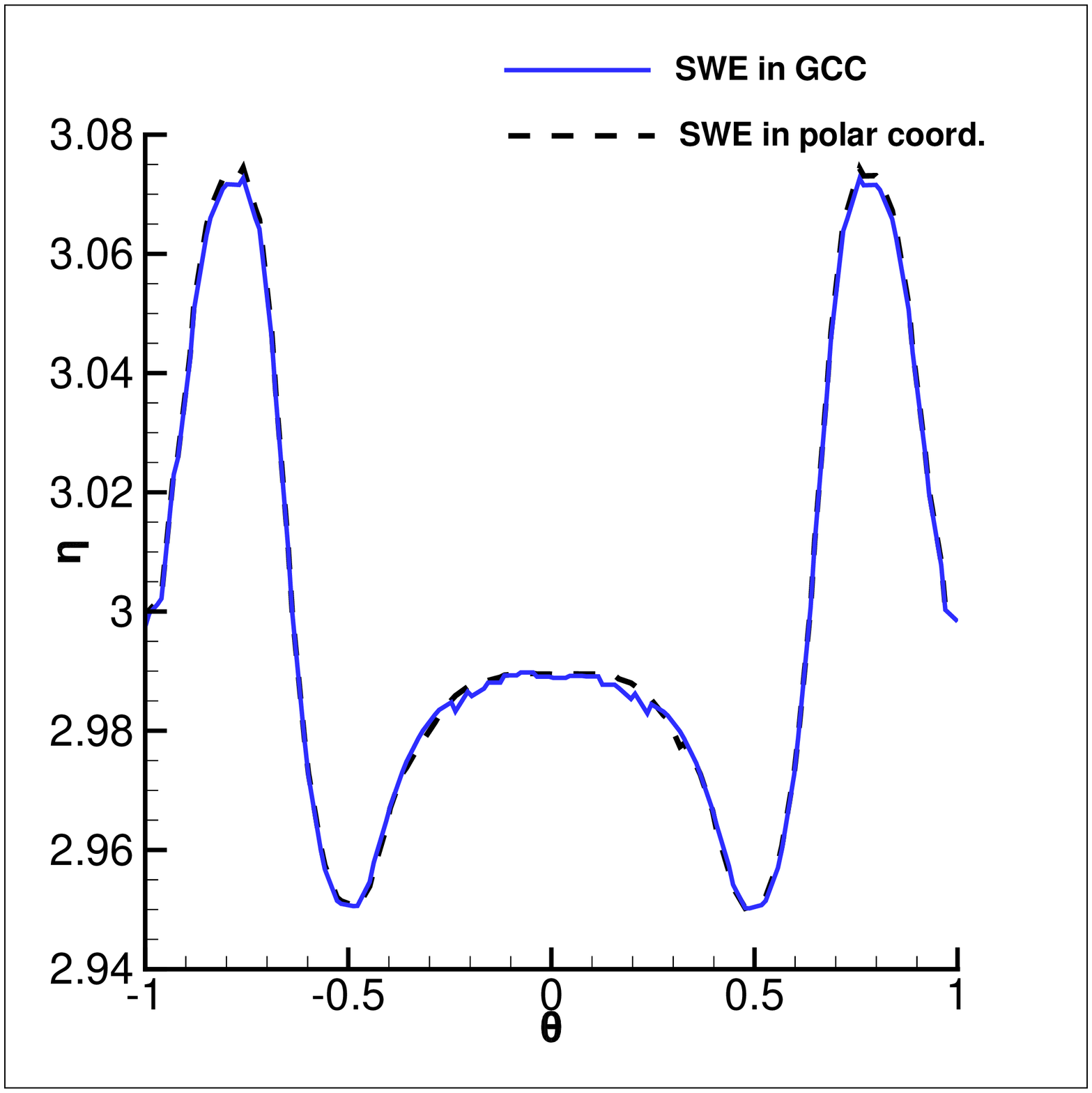}
	}
	\caption{\RIcolor{Comparison of solutions of dambreak-type problem in spherical coordinates from the proposed model~\eqref{eq.hyp} and model~\eqref{eq.castro} presented by Castro \textit{et. al} in~\cite{castro2017well}. In figures (a) and (b), the bottom surface represents the bathymetry and the top surface is the free surface. For both cases, the final time and the averaged mesh size are $T = 0.3$ and $d_N = 5E$-2. In (c), free surface profile of dambreak-type problem at time $T = 0.3$ along the line at $\varphi = 0$ and $\theta \in (-1,1)$. Solid blue line refers to solution from model~\eqref{eq.hyp} with spherical metric. Dashed black line refers to model~\eqref{eq.castro} in~\cite{castro2017well}.}}
	\label{fig.dambreak}
\end{figure}

The last comparison is performed in spherical metric with the model proposed in~\cite{castro2017well}. The hyperbolic formulation of SW equations in spherical coordinates proposed in that work reads
\begin{equation} \label{eq.castro}
	\begin{aligned}
		&\partial_t h_\sigma + \frac{1}{R} \left[ \partial_\theta \left( \frac{M_\theta}{\sigma} \right) + \partial_\varphi M_\varphi \right] = 0, \\
		&\partial_t M_\theta + \frac{1}{R} \partial_\theta \left( \frac{M_\theta^2}{h_\sigma \sigma} \right) + \frac{1}{R} \partial_\varphi \left( \frac{M_\theta M_\varphi}{h_\sigma} \right) + \frac{M_\theta M_\sigma}{R h_\sigma \sigma} \partial_\varphi \sigma + \frac{g h_\sigma}{R \sigma^2} \partial_\theta \eta_\sigma = 0, \\
		&\partial_t M_\varphi + \frac{1}{R} \partial_\theta \left( \frac{M_\varphi M_\theta}{h_\sigma \sigma} \right) + \frac{1}{R} \partial_\varphi \left( \frac{M_\varphi^2}{h_\sigma} \right) - \left( \frac{M_\theta^2}{R h_\sigma \sigma} + \frac{g h_\sigma \eta_\sigma}{R \sigma^2} \right) \partial_\varphi \sigma + \frac{g h_\sigma}{R \sigma} \partial_\varphi \eta_\sigma = 0, \\
		& \partial_t \sigma = 0,
	\end{aligned}
\end{equation}
where the spherical coordinates are the longitudinal-latitudinal angles $(\theta, \varphi)$, $R$ is the radius of the sphere, $\sigma = \cos \varphi$ and the physical quantities read $h_\sigma = h \sigma$ for the fluid depth, $b_\sigma = b \sigma$ for the bathymetry and $M_\Upsilon = m_\Upsilon \sigma$, with $\Upsilon = \theta, \varphi$, for the mass fluxes. 
System~\eqref{eq.castro} and hyperbolic model~\eqref{eq.hyp} in spherical metric proposed in this paper are not comparable because system~\eqref{eq.castro} incorporates curvature and other typical features of bathymetry on the sphere in a slightly different way. For this reason, the spectra of the Jacobian tensors associated with the two systems differ. 
\RIcolor{We report therefore a comparative study of two test cases solved with the two models. 
Concerning the first test case, it is a general Riemann problem. Let us consider the Riemann problem for SW equations on a section of sphere of radius $R = 1$ and defined by angles in $\Omega = [-1,1]^2$ with constant bathymetry $b(\theta, \varphi) \equiv 1$. 
The initial velocity is set to 0 and the initial free surface has a discontinuity along $\theta = 0$. 
In particular, it is defined as $\eta(\theta, \varphi, t = 0) = 2 + 0.5\chi_{\{  \theta > 0\}}(\theta, \varphi)$. Fig.~\ref{fig.SWEPolar} compares the free surface obtained through the proposed hyperbolic model~\eqref{eq.hyp} with metric specialized to be spherical (a) and the free surface by model~\eqref{eq.castro} (b). The colorplots are in good agreement.
The second test case simulates a circular dambreak-type problem. The bathymetry is a sinusoidal hill (e.g. $b(\theta, \varphi) = [0.5 \cos( \pi (\theta^2 + \varphi^2) ) ] \chi_{\{\theta^2 + \varphi^2 \leq 1\}}(\theta, \varphi)$, see bottom surface in Fig. \ref{fig.dambreak}) in the domain $\Omega = (-1.1, 1.1)^2$. At initial time $t = 0$, velocity is zero and the free surface is perturbed by Gaussian bell $\eta(\theta, \varphi) = 10^5 \exp( -1/(0.3^2 - \theta^2 - \varphi^2) ) \chi_{\{ \theta^2 + \varphi^2 \leq 0.3^2 \}}(\theta, \varphi) + 3$. In Fig. \ref{fig.dambreak} (a) and (b), a qualitative comparison for the evolution of the free surface at time $T = 0.3$ from the two models is reported. In Fig. \ref{fig.dambreak} (c) it is depicted the free surface profile along the cut $(\theta, \varphi) \in (-1, 1) \times \{  0 \}$ at time $T = 0.3$. The two curves overlap for both models consistently. The latter test case opens, in the future, the possibility of using the proposed model to simulate, study and analyze geophysical phenomena related to tsunamis on geoids.}

\medskip

This numerically proves the advantage of using the proposed hyperbolic system~\eqref{eq.hyp} that is able to automatically compute the geometrical features of any kind of manifold and its curvature by simply passing as input data the metric tensor components.

\section{Conclusions} \label{sec.concl}
In this paper, we presented a well balanced approach 
for a novel formulation of the shallow water equations on manifolds in general covariant coordinates.
We took as starting point the model by Baldauf in~\cite{baldauf2020discontinuous} and we have proposed a hyperbolic reformulation of the problem. In particular, once the manifold defined by one particular equipotential surface is detected, the induced metric tensor in covariant representation is collected among the conserved variables. 
The obtained model is thus able to automatically compute the curvature of the manifold while the physical quantities are evolved. 
The numerical results prove that the proposed minimally invasive well balanced MUSCL-Hancock approach actually preserves the water at rest solution at machine precision for any nonsingular metric. In addition, it is still second-order accurate for states not of the type water at rest. 
The procedure of locally approximating the fluid depth through the difference between the local reconstructions of the free surface and the bathymetry allows to cancel at discrete level possible numerical problems related to the jumps of the bottom topography. Consequently, if a classical standard non well balanced MUSCL-Hancock approach recovers a numerical solution presenting non physical oscillations in a neighborhood of the discontinuity of the bottom topography, a non oscillating solution (eventually endowed of a smooth free surface) is obtained by the proposed well balanced scheme. 

The idea of involving the curvature of the manifold in the nonconservative components of the formulated hyperbolic system comes from a similar approach for the gravitational field presented in~\cite{gaburro2021well} in the ambit of general relativity. 
In that case indeed the metric plays a fundamental role and it also evolves during time. 
Thus, a first extension of this work will be the development of a minimal invasive well balanced technique also for systems of equations involving general relativity, 
in order to develop a method less expensive than the one presented in~\cite{gaburro2021well}. 
Systems that could benefit of a similar technique could be the GRMHD system~\cite{fambri2018ader,fambri2020discontinuous}, the Einstein field equations~\cite{dumbser2018conformal} and the novel teleparallel formulation of general relativity preliminary introduced in~\cite{Torsion2019,olivares2022new}.

Then, in order to make more effective these techniques in presence of small perturbations of equilibrium solutions, we plan to extend the present methodology to high order of accuracy in an ADER-fashion, following for example~\cite{gaburro2021posteriori,boscheriAFE2022,chiocchetti2021high,han2021dec,gomez2021high}.
And furthermore, since an application of interest, both in shallow water equations and in general relativity is the study of vortical flows for very long simulation times (for energy extraction from water turbines and for the evolution of gas clouds around black holes and neutron stars), the next step of this study will consist in its extension to moving mesh codes specialized in maintaining a high quality of the moving meshes also in presence of strong differential rotation, as those presented in~\cite{springel2010pur,gaburro2017direct,gaburro2020high,gaburro2021unified,bergmann2022ader,moes2021extreme}.

Finally, an important open issue is given by the treatment of the critical points for the metric, namely those points in the computational domain over which the determinant of the metric tensor goes to zero. Different strategies were proposed in the literature: we cite~\cite{arpaia2020well, castro2017well, baldauf2020discontinuous,arpaia2022efficient} among others. 
In our article, we have already seen that it is possible to relax the assumption of smooth equilibrium solutions to be preserved over long timescales. Further analysis, then, can perhaps be undertaken to understand how the system can be adapted when it tries to evolve physical quantities with discontinuous metrics that, globally, describe the immersed manifold without even making the determinant vanish. 
Such approaches would be positioned in a hybrid context between full local approaches, which avoid singularities by referencing each cell to a reference cell, and full global methods, which have the advantage, on one hand, of having the curvature of the manifold computed more easily and automatically and, on the other hand, of modeling also discontinuous bathymetries. Finally, applications to real engineering and geophysical problems will be accounted.

\section*{Acknowledgment}
The authors are members of the CARDAMOM team at the Inria center of the university of Bordeaux.
The authors gratefully acknowledge the support received from the European Union’s Horizon 2020 Research and Innovation Programme under the Marie Skłodowska-Curie Individual Fellowship \textit{SuPerMan}, grant agreement No. 101025563.

\bibliographystyle{plain}
\bibliography{references}

\end{document}